\documentclass[reqno]{amsart}
\usepackage{enumerate}
\usepackage{scrextend}
\usepackage{amsmath,amsxtra,amssymb,latexsym, amscd,amsthm,pb-diagram,fancyhdr,euscript}
\usepackage{amsthm}
\usepackage{latexsym}
\usepackage{amssymb}
\usepackage{amsmath}
\usepackage{indentfirst}
\usepackage{hyperref}
\usepackage{mathrsfs}
\usepackage{array}
\usepackage{euscript}
\usepackage{slashed}
\usepackage[T1]{fontenc}
\usepackage{graphicx}
\usepackage[utf8]{inputenc}
\usepackage[french,english]{babel}
\usepackage{lipsum} 
\usepackage{multicol}
\hypersetup{
   colorlinks,
    citecolor=blue,
    filecolor=black,
    linkcolor=blue,
    urlcolor=magenta
}
\hypersetup{linktocpage}
\changefontsizes[16pt]{12pt}
\newtheorem{theorem}{Theorem}[section]
\newtheorem{lemma}[theorem]{Lemma}

\theoremstyle{definition}
\newtheorem{definition}[theorem]{Definition}

\newtheorem{remark}[theorem]{Remark}

\newcommand{\norm}[1]{\left\Vert#1\right\Vert}

\numberwithin{equation}{section}
\usepackage{amsfonts}
\usepackage{amsmath}
\usepackage{amssymb}

\begin{document}
\font\nho=cmr10
\def\dive{\mathrm{div}}
\def\cal{\mathcal}
\def\L{\cal L}

\def \ud{\underline }
\def\id{{\indent }}
\def\f{\frac}
\def\non{{\noindent}}
 \def\le{\leqslant} 
 \def\leq{\leqslant}
 \def\geq{\geqslant} 
\def\rar{\rightarrow}
\def\Rar{\Rightarrow}
\def\ti{\times}
\def\i{\mathbb I}
\def\j{\mathbb J}
\def\si{\sigma}
\def\Ga{\Gamma}
\def\ga{\gamma}
\def\ld{{\lambda}}
\def\Si{\Psi}
\def\f{\mathbf F}
\def\r{\hro{R}}
\def\e{\cal{E}}
\def\B{\cal B}
\def\A{\mathcal{A}}
\def\p{\mathbb P}

\def\tet{\theta}
\def\Tet{\Theta}
\def\hro{\mathbb}
\def\ho{\mathcal}
\def\P{\ho P}
\def\E{\mathcal{E}}
\def\n{\mathbb{N}}
\def\M{\mathbb{M}}
\def\dMu{\mathbf{U}}
\def\dMcs{\mathbf{C}}
\def\dMcu{\mathbf{C^u}}
\def\vk{\vskip 0.2cm}
\def\td{\Leftrightarrow}
\def\df{\frac}
\def\Wei{\mathrm{We}}
\def\Rey{\mathrm{Re}}
\def\s{\mathbb S}
\def\l{\mathcal{L}}
\def\C+{C_+([t_0,\infty))}
\def\o{\cal O}

\title[AAP- Mild Solutions of INSE]{On asymptotically almost periodic solutions to the Navier-Stokes equations on hyperbolic manifolds}

\author[P.T. Xuan]{Pham Truong Xuan}
\address{Pham Truong Xuan \hfill\break
Corresponding author, Faculty of Pedagogy, VNU University of Education, Vietnam National University, 144 Xuan Thuy, Cau Giay, Hanoi, Viet Nam} 
\email{phamtruongxuan.k5@gmail.com or ptxuan@vnu.edu.vn} 

\author[N.T. Van]{Nguyen Thi Van}
\address{Nguyen Thi Van \hfill\break
Faculty of Information Technology, Department of Mathematics, Thuyloi university \hfill\break
Khoa Cong nghe Thong tin, Bo mon Toan, Dai hoc Thuy loi, 175 Tay Son, Dong Da, Ha Noi, Viet Nam}
\email{van@tlu.edu.vn}

\begin{abstract}  
In this paper, we study the forward asymptotically almost periodic (AAP-) mild solutions of Navier-Stokes equations on the real hyperbolic manifold $\mathcal{M}=\mathbb{H}^d(\mathbb{R})$ with dimension $d \geq 2$. Using the dispersive and smoothing estimates for the Stokes equation, we invoke the Massera-type principle to prove the existence and uniqueness of the AAP- mild solution for the inhomogeneous Stokes equations in $L^p(\Gamma(T\mathcal{M})))$ space with $1<p\leq d$. Next, we establish the existence and uniqueness of the small AAP- mild solutions of the Navier-Stokes equations by using the fixed point argument, and the results of inhomogeneous Stokes equations. The asymptotic behaviour (exponential decay and stability) of these small solutions are also related. This work, together with our recent work [P.T. Xuan, N.T. Van and B. Quoc, {\it On Asymptotically Almost Periodic Solution of Parabolic Equations on real hyperbolic Manifolds},  J. Math. Anal. Appl., Vol.
517, Iss. 1 (2023), pages 1-19], provide a full existence and asymptotic behaviour of AAP- mild solutions of Navier-Stokes equations in $L^p(\Gamma(T\mathcal{M})))$ spaces for all $p>1$.
\end{abstract}

\subjclass[2010]{Primary 35Q30, 35B35, 43A60, 34K25; Secondary 32Q45, 58J35}

\keywords{Navier-Stokes Equations, hyperbolic manifold, Deformation tensor, Asymptotically Almost periodic functions (resp. solutions), exponential decay (stability)}

\maketitle

\tableofcontents

\section{Introduction and preliminaries}
\subsection{Introduction}
The incompressible Navier-Stokes equations (INSE) have been widely studied on the hyperbolic or non-compact Einstein manifolds with negative Ricci curvatures since 1970's. The original paper of Ebin and Marsden \cite{EbiMa} established the formula of Navier-Stokes equation on these manifolds. In this work, the authors gave a generalized Laplace operator on vector fields on Einstein manifolds by using the deformation tensor. Then, Ebin and Marsden's Laplace operator has been studied extensively in the works of Czubak and Chan \cite{Cz1,Cz2,Cha},  Pierfelice \cite{Pi}, Nguyen et al. \cite{HuyXuan2020,HuyXuan2021,HuyXuan2021'} and many other authors (see for example \cite{Ba2020,Fa2018,Fa2020,Li,Zha}). In particular, by using Ebin-Marsden's Laplace operator, Czubak and Chan \cite{Cz1,Cz2} (see also \cite{Khe2012,Li}) considered the problem of non-uniqueness of the weak Leray-Holf solutions of INSE on the real hyperbolic manifold, or some other non-compact manifolds. Pierfelice \cite{Pi} established the dispersive and smoothing estimates, then proved the existence and uniqueness of bounded mild solutions of INSE on non-compact Einstein and generalized non-compact manifolds with negative Ricci curvatures. 

The study of periodic mild solution and its generalizations of INSE is an important direction, and it has received a lot of attention from mathematicians. On Euclidean space $\mathbb{R}^n$ or its unbounded domains, the various types of bounded mild solution of INSE have been studied in some recent works \cite{GHN,HiHuSe,Huy2014,Huy2018,HuyHaSacXuan2,HuyHaSacXuan2',HuyHaSacXuan3}. In these works, the authors used the $L^p-L^q$- dispersive and smoothing estimates of Stokes semigroups to invoke the Massera-type principle, then prove the existence, uniqueness and polynomial stability of the the periodic, almost periodic, almost automorphic and some other mild solutions of INSE in certain interpolation spaces.
However, up to our knowledge, for the asymptotically almost periodic (AAP-) mild solutions of INSE there is one work done by Farwig and Tanuichi (see \cite{FaTa2013}) which gives the global uniqueness of the small backward AAP- mild solution in unbounded domains in $\mathbb{R}^n$. There are no results on the existence, uniqueness and stability of the forward AAP-mild solution of INSE in $\r^n$ so far. The challenges come from the fact that when we consider the Stokes equations in $\mathbb{R}^n$ or its unbounded domains, the semigroups associated with equations are polynomial stable, hence the solution operator (see its definition in Section \ref{S3}) does not preserve the forward asymptotic property of the AAP- functions. Therefore, the existence, uniqueness and stability of the forward AAP- mild solution of INSE on $\mathbb{R}^n$, or its unbounded domains are remain an open problem.

In a recent work \cite{XVQ2023}, we have answered partly this problem in the case of real hyperbolic manifold $\mathcal{M}=\mathbb{H}^d(\mathbb{R})$ with $d\geqslant 2$. In particular, we have studied the generalized parabolic evolution equations (which consist of Navier-Stokes equations) on $\mathcal{M}$, and established the existence, uniqueness and exponential stability of forward AAP- mild solutions in the space $L^p(\Gamma(T\mathcal{M}))$ for all $p>d$. 

In other related works, the existence, uniqueness and exponential stability of the periodic and almost periodic mild solutions of INSE on the non-compact Einstein or generalized non-compact manifolds are established by Nguyen and his collaborations in \cite{HuyXuan2020, HuyXuan2021, HuyXuan2021'}. The authors combined the $L^p-L^q$- dispersive and smoothing estimates with the ergodic method to invoke the Massera-type principle, then construct the initial data and prove the existence, uniqueness of the periodic and almost periodic mild solutions of the vectorial heat or Stokes equations. By using these existence results, fixed point arguments and the cone inequality, the authors proved the existence, uniqueness and exponential stability of such solutions of INSE.
  
In this paper, we extend results obtained in \cite{XVQ2023} to study the existence and uniqueness of the forward AAP- mild solution of INSE in the space $L^p(\mathcal{M})$ for $1<p\leqslant d$. We would like to note that, on hyperbolic space $\mathcal{M}=\mathbb{H}^d(\mathbb{R})$, the semigroup generated by Stokes equation is exponential stable (see Lemma \ref{estimates} or more detail in \cite[Theorem 4.1 and Corollary 4.3]{Pi}). This leads to the fact that, the solution operator preserve the forward asymptotic property of the functions. Therefore, we will use the $L^p-L^q$- dispersive and smoothing estimates of Pierfelice to invoke the Massera-type principle and prove the existence of the forward AAP- mild solution for the inhomogeneous Stokes equation (see \cite{HiHuSe,Huy2014,Huy2018,HuyHaSacXuan3}, for the detailed method for INSE on unbounded domains in $\mathbb{R}^n$ and \cite{HuyXuan2020,HuyXuan2021'} for the case of Einstein or generalized non-compact Riemannian manifolds). In particular, we prove an argument that, if the external force of INSE is a forward AAP- function, then the bounded mild solution of the corresponding inhomogeneous Stokes equation is also a forward AAP- function. Then, we use the fixed point argument to prove the existence and uniqueness of AAP- mild solution of INSE. The exponential decay will be obtained by using the Gronwall's inequality, and then the exponential stability will be obtained as a direct consequence of this decay. 

Our results and the ones obtained in \cite{XVQ2023} answer completely the questions about the existence, uniqueness and stability of the forward AAP- mild solution of INSE in $L^p(\mathcal{M})$ spaces for all $p>1$. This paper is the continuation of the recent works to study INSE, as well as parabolic evolution equations on the non-compact Riemannian manifolds with negative Ricci curvatures \cite{HuyXuan2020,HuyXuan2021,HuyXuan2021',XVQ2023}.

Our paper is organised as follows: In Section \ref{S2}, we recall the formula of INSE on the real hyperbolic manifold; Sections \ref{S3}, we prove the existence and uniqueness of the forward AAP- mild solutions for the inhomogeneous Stokes equation; Section \ref{S4}, we establish the existence, uniqueness and exponential decays (stability) of the small forward AAP- mild solutions of INSE; Appendix \ref{A} gives some detailed calucations which are useful for previous sections.

\subsection{Preliminaries}
Before starting the main sections, we recall some basic definitions and give some notations. Throughout this article, we denote 
$$C_b(\r, X):=\{f:\r \to X \mid f\hbox{ is continuous on $\r$ and }\sup_{t\in\r}\|f(t)\|_X<\infty\},$$
in which, $X$ is  a  Banach space with norm $\|f\|_{\infty, X}=\|f\|_{C_b(\r, X)}:=\sup_{t\in\r}\|f(t)\|_X$.
We recall the definitions of almost periodic and asymptotically almost periodic functions (for details see \cite{Che,Dia}).
\begin{definition}
A function  $h \in C_b(\r, X )$ is called almost periodic function if for each $ \epsilon  > 0$, there exists $l_{\epsilon}>0 $ such that every interval of length $l_{\epsilon}$ contains at least a number $T $ with the following property
\begin{equation}
 \sup_{t \in \r } \| h(t+T)  - h(t) \| < \epsilon.
\end{equation}
The collection of all almost periodic functions $h:\r \to X $ will be denoted by $AP(\r,X)$ which is a Banach space endowed with the norm $\|h\|_{ AP(\r,X)}=\sup_{t\in\r}\|h(t)\|_X.$
\end{definition}
To introduce the forward asymptotically almost periodic functions, we need the space  $C_0 (\r_+,X)$, that is the collection of all forward asymptotic and continuous functions $\varphi: \r_+ \to X$ satisfying
$$\lim_{t \to +\infty } \| \varphi(t) \|=0.$$
Clearly, $C_0 (\r_+,X)$  is a Banach space endowed with the norm $\|\varphi\|_{C_0 (\r_+,X)}=\sup_{t\in\r_+}\|\varphi(t)\|_X$.
\begin{definition} 
A continuous function  $f \in C(\r_+, X )$  is said to be forward asymptotically almost periodic if there exist  $h \in AP(\r,X)$ and $ \varphi\in C_0(\r_+,X)$ such that
\begin{equation}
f(t) = h(t) + \varphi(t).
\end{equation}
We denote $AAP(\r_+, X):= \{f:\r_+ \to X \mid f\hbox{ is asymptotically almost periodic on $\r_+$}\}$. Note that $AAP(\r_+,X)$ is a Banach space with the norm defined by $\|f\|_{ AAP(\r_+,X)}=\|h\|_{ AP(\r, X)}+\|\varphi\|_{ C_0(\r_+,X)}$.
\end{definition}
We notice that, if $\varphi: \r_- \to X$ is a backward asymptotic and continuous function, i.e
$$\lim_{t \to -\infty } \| \varphi(t) \|=0$$
and $h\in AP(\r,X)$, then $f(t) = h(t) + \varphi(t)$ is a backward asymptotcally almost periodic function (see \cite{FaTa2013}).\\
{\bf Notations.}\\
$\bullet$ We denote the Levi-Civita connection by $\nabla$; $\Gamma(T\mathcal{M})$ for the set of all vector fields on $\mathcal{M}$, and $\Gamma(T\mathcal{M}\otimes T\mathcal{M})$ for the set of all second order tensor fields on $\mathcal{M}$.\\
$\bullet$ We shall use the notation $\lesssim$ in the sense that $f\lesssim g$ if and only if there exists a positive constant $C$ which is independent of $f$ and $g$ such that $f\leq C g$.\\

\section{Setting of the incompressible Navier-Stokes Equation on the real hyperbolic Manifold}\label{S2}
Let  $(\mathcal{M} =\mathbb{H}^d(\mathbb{R}),g)$  be a real hyperbolic manifold of dimension $d\geq 2$ which is realized as the upper sheet 
$$x_0^2-x_1^2-x_2^2...-x_d^2 = 1 \,  \,( x_0\geq 1),$$
of hyperboloid in $\mathbb{R}^{d+1}$, equipped with the Riemannian metric 
$$g = -dx_0^2 + dx_1^2 + ... + dx_d^2.$$
In geodesic polar coordinates, the hyperbolic manifold is 
$$\mathbb{H}^d(\mathbb{R}) = \left\{ (\cosh \tau, \omega \sinh \tau), \, \tau\geq 0, \omega \in \mathbb{S}^{d-1}  \right\}$$
with the metric 
$$g = d\tau^2+(\sinh\tau)^2d\omega^2,$$
where  $d\omega^2$ is the canonical metric on the sphere $\mathbb{S}^{d-1}$. 
 
A remarkable property on $\mathcal{M}$ is that: $\mathrm{Ric}_{ij}=-(d-1)g_{ij}$, where $\mathrm{Ric}_{ij}$ is the component of Ricci curvature tensor.

We follow \cite{HuyXuan2020,Pi} to express the Navier-Stokes equation on the real hyperbolic manifold. Denote the Levi-Civita connection by $\nabla$ and the set of vector fields on $\mathcal{M}$ by $\Gamma(T\mathcal{M})$. The imcompressible Navier-Stokes equation on $\mathcal{M}$ with the initial data $u|_{t=0}=u(0)$ are described by the following system:
\begin{align}\label{Cauchy}
\begin{cases}
\partial_t u + \nabla_u u + \nabla p = Lu + \mathrm{div}F \cr
\mathrm{div} u = 0\cr
u|_{t=0} = u(0) \in \Gamma(T\mathcal{M}),
\end{cases}
\end{align}
where $u=u(x,t)$ is considered as a vector field on $\mathcal{M}$, i.e., 
$u(\cdot,t)\in \Gamma(T\mathcal{M})$ and $p=p(x,t)$ is the pressure field, $L$ is the stress tensor, $\dive F$ is the external force, where $F(\cdot,t)\in \Gamma(T\mathcal{M}\otimes T\mathcal{M})$ is a second order tensor field on $\mathcal{M}$.

Since $\mathrm{div}u=0$, it follows that $\nabla_uu = \mathrm{div} (u\otimes u)$. Besides, the vectorial Laplacian $L$ is defined by using the deformation tensor 
$$Lu := \frac{1}{2}\mathrm{div}(\nabla u + \nabla u^t)^{\sharp},$$
where $\omega^{\sharp}$ is the vector field associated with the $1$-form $\omega$ by $g(\omega^{\sharp},Y) = \omega(Y)$ for all $Y \in \Gamma(T\mathcal{M})$.
Since $\mathrm{div} u=0$ , $L$ can be expressed by
$$Lu = \overrightarrow{\Delta}u + R(u),$$
where $\overrightarrow{\Delta}u =- \nabla^*\nabla u= \mathrm{Tr}_g(\nabla^2u)$ is the Bochner-Laplacian
and $R(u)=(\mathrm{Ric}(u,\cdot))^{\sharp}$ is the Ricci operator. Since $\mathrm{Ric}(u,\cdot)=-(d-1)g(u,\cdot)$, we have $R(u)=-(d-1)u$. Therefore, we obtain the equivalent system of \eqref{Cauchy} as follows
\begin{align}\label{Cauchy'}
\begin{cases}
\partial_t u + \mathrm{div}(u\otimes u) + \nabla p = \overrightarrow{\Delta}u - (d-1)u + \mathrm{div}F \cr
\mathrm{div} u = 0\cr
u|_{t=0} = u(0) \in \Gamma(T\mathcal{M}).
\end{cases}
\end{align}

By the Kodaira-Hodge decomposition, an $L^2-$ form can be decomposed  on $\mathcal{M}$ as 
$$L^2(\Gamma(T^*(\mathcal{M}))) = \overline{\mathrm{Image}\, d} \oplus \overline{\mathrm{Image}\, d^*} \oplus \mathcal{H}^1(\mathcal{M}) \, ,$$
where $\mathcal{H}^1(\mathcal{M})$ is the space of $L^2$ harmonic $1-$forms. Taking the divergence of equations \eqref{Cauchy'}, and noting that $\mathrm{div}(\overrightarrow{\Delta}u) = \mathrm{div}(-(d-1)u) = 0$ if $\mathrm{div}u=0$, we get
$$\Delta_g p + \mathrm{div}[\nabla_uu] = 0,$$
where $\Delta_g$ is the Laplace-Beltrami operator on $\mathcal{M}$.

\def\dive{\operatorname{div}}
We need to choose a solution in $L^p$ of this elliptic equation. Since the spectral of $\Delta_g$ on the hyperbolic manifold $\mathcal{M} = \mathbb{H}^d(\mathbb{R})$ is $\left[\frac{(d-1)^2}{4}, \infty \right)$ which does not contain $0$, we have that
 $\Delta_g \,: \, W^{2,r} \rightarrow L^r$ is an isomorphism for $2\leq r<\infty$ and 
$$\mathrm{grad}p = \mathrm{grad}(-\Delta_g)^{-1}\mathrm{div}[\nabla_uu].$$
By setting the operator 
$\mathbb{P}:= I + \mathrm{grad}(-\Delta_g)^{-1} \mathrm{div}$, we can get rid of the 
pressure term $p$ and then obtain from \eqref{Cauchy'} the following system 
\begin{align}\label{DivNavierStokes}
\begin{cases}
\partial_t u - (\overrightarrow{\Delta}u -(d-1)u) = -\mathbb{P}\dive (u\otimes u) + \mathbb{P}\dive F\\
\dive u = 0\cr
u|_{t=0} = u(0) \in \Gamma(T\mathcal{M}).
\end{cases}
\end{align}

\section{Forward asymptotically almost periodic mild solution of the Stokes equation}\label{S3}
We consider the inhomogeneous Stokes equation
\begin{align}\label{CauchyHeat}
\begin{cases}
\partial_t u = \overrightarrow{\Delta}u - (d-1)u - \p \dive(v\otimes v) + \p\dive F\cr
\dive u=0,
\end{cases}
\end{align}
for a given vector field $v(t,\cdot)$ in $\Gamma(T\mathcal{M})$ and a given second order tensor field $F(t,\cdot)$ in $\Gamma(T\mathcal{M}\otimes T\mathcal{M})$.

Setting $\mathcal{A}u  = -(\overrightarrow{\Delta}u - (d-1)u)$, and denote 
$e^{-t\A}$ the semigroup associated with the homogeneous Cauchy problem of the vectorial heat equation
$$\partial_t u = -\mathcal{A}u.$$

We now recall the results on the dispersive and smoothing estimates (see \cite{Pi}) for the semigroup $e^{-t\A}$:
\begin{lemma}{\rm(\cite[Theorem 4.1, Corollary 4.3]{Pi}):}\label{estimates}
\mbox{}\\
\begin{enumerate}[(i)]
\item For $t>0$, and $p$, $q$ satisfying $1\leq p \leq q \leq \infty$, 
the following dispersive estimates hold 
\begin{equation}\label{dispersive}
\left\| e^{-t\mathcal{A}} u(0)\right\|_{L^q} \leq [h_d(t)]^{\frac{1}{p}-\frac{1}{q}}e^{-t(d-1 + \gamma_{p,q})}\left\| u(0) \right\|_{L^p}\hbox{ for all }u(0) \in L^p(\Gamma(T\mathcal{M}))
\end{equation}
 where $h_d(t): = C\max\left( \frac{1}{t^{d/2}},1 \right)$, 
   $\gamma_{p,q}:=\frac{\delta_d}{2}\left[ \left(\frac{1}{p} - \frac{1}{q} \right) + \frac{8}{q}\left( 1 - \frac{1}{p} \right) \right]$ and $\delta_d$ are a positive constants depending only on $d$.  
\item For $p$ and $q$ satisfying $1<p\leq q <\infty$, we obtain that, for all $t>0$:
\begin{equation}
\left\| \nabla e^{-t\mathcal{A}}u(0) \right\|_{L^q} \leq [h_d(t)]^{\frac{1}{p}-\frac{1}{q}+\frac{1}{d}}e^{-t\left(d-1  + \frac{\gamma_{q,q}+\gamma_{p,q}}{2} \right)} \left\| u(0) \right\|_{L^p},
\end{equation}
\begin{equation}
\left\| e^{-t\mathcal{A}}\nabla^* T_0 \right\|_{L^q} \leq [h_d(t)]^{\frac{1}{p}-\frac{1}{q}+\frac{1}{d}}e^{-t\left(d-1  + \frac{\gamma_{q,q}+\gamma_{p,q}}{2} \right)} \left\| T_0 \right\|_{L^p},
\end{equation}
for all vector fiels $u(0) \in L^p(\Gamma(T\mathcal{M}))$ and all 
 tensor fields $T_0 \in L^p(\Gamma(T\mathcal{M} \otimes T^*\mathcal{M}))$.
\item  As a consequence of (ii), for all $t>0$, we have 
\begin{equation}
\left\| e^{-t\mathcal{A}}\mathrm{div}T^{\sharp}_0 \right\|_{L^q} \leq [h_d(t)]^{\frac{1}{p}-\frac{1}{q}+\frac{1}{d}}e^{-t\left(d-1 + \frac{\gamma_{q,q}+\gamma_{p,q}}{2} \right)} \left\| T_0^{\sharp}\right\|_{L^p},
\end{equation}
for all tensor fields $T^{\sharp}_0 \in L^p(\Gamma(T\mathcal{M}\otimes T\mathcal{M}))$.
\end{enumerate}
\end{lemma}
The mild solution of Equation \eqref{CauchyHeat} on the whole time-line is satisfied the following integral equation (see \cite{KoNa}):
\begin{equation}\label{mild:linear1}
u(t) = \int_{-\infty}^te^{-(t-\tau)\mathcal{A}} \mathbb{P} \dive (-v(\tau)\otimes v(\tau)+ F (\tau)) d\tau \hbox{         } (t\in \mathbb{R}).
\end{equation}
The mild solution of Equation \eqref{CauchyHeat} on the half time-line with the initial data $u|_{t=0}=u(0,\cdot)=u(0) \in \Gamma(T\mathcal{M})$ is satisfied the following integral equation
\begin{equation}\label{mild:linear}
u(t) = e^{-t\mathcal{A}}u(0) + \int_{0}^te^{-(t-\tau)\mathcal{A}} \mathbb{P} \dive (-v(\tau)\otimes v(\tau)+ F (\tau)) d\tau \hbox{          } (t\in \mathbb{R}_+).
\end{equation}

Let $0<\delta<1$ and $0<\sigma<{\beta}= d-1+ \frac{\gamma_{d/\delta,d/\delta}+\gamma_{d/(2\delta),d/\delta}}{2}$, we consider the existence and uniqueness of boundedness (in time) of mild solution on the half time-line axis to Equation \eqref{DivNavierStokes} on the following Banach space
\begin{eqnarray*}
\mathcal{X}&=&\left\{v\in C_b(\r_+, L^p(\Gamma(T\mathcal{M}))\cap L^d(\Gamma(T\mathcal{M})) \cap L^{d/\delta}(\Gamma(T\mathcal{M}))):\right.\cr
&&\hspace{3cm}\left. \|v(t)\|_{p} + \|v(t)\|_{d} + [h_d(t)]^{-\frac{1-\delta}{d}}e^{\sigma t}\norm{v(t)}_{\frac{d}{\delta}} \le \rho \right\},
\end{eqnarray*}
equipped with the norm 
\begin{equation}\label{space1}
\norm{v}_{\mathcal{X}} = \sup_{t>0} \left( \|v(t)\|_{p} + \| v(t) \|_{d} + [h_d(t)]^{-\frac{1-\delta}{d}}e^{\sigma t}\norm{v(t)}_{\frac{d}{\delta}}\right).
\end{equation}
We also consider the existence and uniqueness of boundedness (in time) of mild solution on the full time-line axis to Equation \eqref{CauchyHeat} on 
\begin{eqnarray*}
\mathbb{X}&=&\left\{v\in C_b(\r, L^p(\Gamma(T\mathcal{M}))\cap L^d(\Gamma(T\mathcal{M})) \cap L^{d/\delta}(\Gamma(T\mathcal{M}))):\right.\cr
&&\hspace{5cm}\left. \|v(t)\|_{p} + \|v(t)\|_{d} + \lambda(t)\norm{v(t)}_{\frac{d}{\delta}} \le \rho \right\},
\end{eqnarray*}
where the time-weighted function $\lambda(t)$ is given by
\begin{align}\label{lambda}
\lambda(t) = \begin{cases}
\displaystyle  [h_d(t)]^{-\frac{1-\delta}{d}}e^{\sigma t}&\hbox{   if  } t\geqslant 0 \cr
[h_d(|t|)]^{-\frac{1-\delta'}{d}} & \hbox{   if  } t<0,
\end{cases}
\end{align}
with a given constant $\delta'$ satisfying $0<\delta<\dfrac{1+\delta}{2}<\delta'<1$. We define the norm on the space $\mathbb{X}$ by 
\begin{equation}\label{space2}
\norm{v}_{\mathbb{X}} = \sup_{t\in \mathbb{R}} \left( \|v(t)\|_{p} + \| v(t) \|_{d} + \lambda(t)\norm{v(t)}_{\frac{d}{\delta}}\right).
\end{equation}

We state and prove the existence and uniqueness of the inhomogeneous Stokes equation \eqref{CauchyHeat} on $d$-dimensional real hyperbolic manifold $(\mathcal{M},g)\, (d\geqslant 2)$ in the following lemma:
\begin{lemma}\label{Thm:linear}
Let $(\mathcal{M},g)$ be a $d$-dimensional real hyperbolic manifold with $d \geqslant 2$.
Let $0<\delta<1$ and $0<\sigma<{\beta}= d-1+ \frac{\gamma_{d/\delta,d/\delta}+\gamma_{d/(2\delta),d/\delta}}{2}$. The following assertions hold: 
\begin{itemize}
\item[i)] Suppose that $u(0)\in L^p(\Gamma(T\mathcal{M}))\cap L^d(\Gamma(T\mathcal{M})),\, v\in \mathcal{X}$ and $F$  belongs  to $ C_b(\r_+, (L^{\frac{dp}{d+\delta p}}\cap L^{\frac{d}{1+\delta}}\cap L^{\frac{d}{2\delta}})(\Gamma(T\mathcal{M}\otimes T\mathcal{M})))$ for $1<p\leq d$ and
$$P=\sup\limits_{t >0} \left( \| F(t)\|_{\frac{dp}{d+\delta p}} + \| F(t) \|_{\frac{d}{1+\delta}} + [h_d(t)]^{-\frac{1-\delta}{d}}e^{\sigma t}\|F(t)\|_{\frac{d}{2\delta}} \right)<+\infty.$$
Then, Problem \eqref{CauchyHeat} with the initial value $u(0)$ has one and only one mild solution $u\in \mathcal{X}$ given by the formula \eqref{mild:linear}.
Moreover,
\begin{equation}\label{CoreEstimates}
\norm{u}_{\mathcal{X}}\le \norm{u(0)}_{L^p\cap L^d} +M\norm{v}^2_{\mathcal{X}}+ NP,
\end{equation}
where constants $M$ and $N$ are independent of $u(0), \,v$ and $F$.
\item[ii)] Suppose that $v \in \mathbb{X}$ and $F$  belongs  to $ C_b(\r, (L^{\frac{dp}{d+\delta p}}\cap L^{\frac{d}{1+\delta}}\cap L^{\frac{d}{2\delta}})(\Gamma(T\mathcal{M}\otimes T\mathcal{M})))$ for $1<p\leq d$ and
$$\tilde{P}=\sup\limits_{t \in \mathbb{R}} \left( \| F(t)\|_{\frac{dp}{d+\delta p}} + \| F(t) \|_{\frac{d}{1+\delta}} + \lambda(t)\|F(t)\|_{\frac{d}{2\delta}} \right)<+\infty.$$
Then, Problem \eqref{CauchyHeat} has one and only one mild solution $u\in  
\mathbb{X}$ given by the formula \eqref{mild:linear1}.
Moreover, 
\begin{equation}\label{rgl1}
\norm{u}_{\mathbb{X}}\le \tilde{M}\norm{v}^2_{\mathbb{X}} + \tilde{N}\tilde{P},
\end{equation}
where constants $\tilde{M}$ and $\tilde{N}$ are independent of $v$ and $F$.
\end{itemize} 
Here, we use the notation $L^{r}(Z)\cap L^{s}(Z) \cap L^{z}(Z)= (L^r\cap L^s\cap L^z)(Z)$, where $Z$ is $\Gamma(T\mathcal{M})$ or $\Gamma(T\mathcal{M}\otimes T\mathcal{M})$.
\end{lemma}
\begin{proof}
\def\A{\mathcal{A}}
\noindent
\item[i)] It is clear that, if $u(t)$ satisfies the integral equation \eqref{mild:linear}, then $u(t)$ satisfies Equation \eqref{CauchyHeat}. We need only to estimate the norm $\|u(t)\|_{\mathcal{X}}$ for each $t\in \r_+$. For the sake of convenience, we denote $\|\cdot\|_{r}:=\|\cdot\|_{L^r(Z)}$ and $\norm{\cdot}_{C_b(\r_+, L^r(Z))}=\norm{\cdot}_{\infty,r}$ for $r>0$, and $Z$ is $\Gamma(T\mathcal{M})$ or $\Gamma(T\mathcal{M}\otimes T\mathcal{M})$.
 
Notice that, on the hyperbolic space, we have the $L^p-$boundedness of Riesz transform on the real hyperbolic manifold (see \cite{Loho}). Therefore, the operator $\mathbb{P}$ is bounded. Since the Ricci curvature is constant, operator $\mathbb{P}$ is commuted with $\overrightarrow{\Delta}$ and $e^{-t\A}$. 

First, we estimate $\norm{u(t)}_p$. By using assertions $i)$ and $iii)$ of Lemma \ref{estimates}, we have
\begin{eqnarray}\label{Part1}
\|u(t)\|_{p} &\le& \| e^{-t\A}u(0)\|_{p} + \int_0^t \|e^{-(t-\tau)\A}\dive (-v(\tau)\otimes v(\tau)+ F(\tau))\|_{p}d\tau\cr 
&\le& \| u(0)\|_{p} + \int_0^t [h_d(t-\tau)]^{\frac{\delta}{d}+\frac{1}{d}}e^{-\left(d-1  + \frac{\gamma_{p,p} + \gamma_{dp/(1+\delta p),p}}{2} \right)(t-\tau)} \norm{v(\tau)\otimes v(\tau)}_{\frac{dp}{d+\delta p}} d\tau\cr
&&+ \int_0^t [h_d(t-\tau)]^{\frac{\delta}{d}+\frac{1}{d}}e^{-\left(d-1  + \frac{\gamma_{p,p} + \gamma_{dp/(1+\delta p),p}}{2} \right)(t-\tau)} \left\| F(\tau) \right\|_{\frac{dp}{d+\delta p}} d\tau\cr
&\le& \| u(0)\|_{p} + \int_0^t [h_d(t-\tau)]^{\frac{1+\delta}{d}}e^{-\beta_1(t-\tau)}\norm{v(\tau)}_{\frac{d}{\delta}}\norm{v(\tau)}_p d\tau\cr
&&+\int_0^t C^{ \frac{\delta}{d}+\frac{1}{d}} \left((t-\tau)^{-\frac{1+\delta}{2}}+1\right) e^{-\beta_1 (t-\tau)} d\tau \|F\|_{\infty,\frac{dp}{d+\delta p}}\cr
&\le& \| u(0)\|_{p} + \norm{v}^2_{\mathcal{X}}\int_0^t [h_d(t-\tau)]^{\frac{1+\delta}{d}}[h_d(\tau)]^{\frac{1-\delta}{d}}e^{-\beta_1(t-\tau)}e^{-\sigma\tau}d\tau\cr
&&+\int_0^t C^{ \frac{1+\delta}{d}} \left((t-\tau)^{-\frac{1+\delta}{2}}+1\right) e^{-\beta_1 (t-\tau)} d\tau \|F\|_{\infty,\frac{dp}{d+\delta p}}\cr
&\le& \|u(0)\|_{p} + M_1\norm{v}^2_{\mathcal{X}} + \tilde{C}\left(\beta^{\theta-1}_1{\mathbf{\Gamma}}(1-\theta) + \frac{1}{\beta_1} \right)P,
\end{eqnarray}
where $\beta_1 = d-1  + \frac{\gamma_{p,p} + \gamma_{dp/(1+\delta p),p}}{2}$, $\tilde{C}=C^{\frac{1+\delta}{d}}>0$, $0<\theta=\dfrac{1+\delta}{2}<1$, function ${\mathbf{\Gamma}}$ is the gamma function, and the integral 
$$\int_0^t [h_d(t-\tau)]^{\frac{1+\delta}{d}}[h_d(\tau)]^{\frac{1-\delta}{d}}e^{-\beta_1(t-\tau)}e^{-\sigma\tau}d\tau \leqslant M_1 <+\infty$$
converges by the mean of gamma and beta functions (see Appendix).

By the same way, we estimate $\norm{u(t)}_d$. In details, 
\begin{eqnarray}\label{Part2}
\|u(t)\|_{d} &\le& \| e^{-t\A}u(0)\|_{d} + \int_0^t \|e^{-(t-\tau)\A}\dive (-v(\tau)\otimes v(\tau)+ F(\tau))\|_{d}d\tau\cr 
&\le& \| u(0)\|_{d} + \int_0^t [h_d(t-\tau)]^{\frac{1+\delta}{d}}e^{-\left(d-1  + \frac{\gamma_{d,d} + \gamma_{d/(1+\delta),d}}{2} \right)(t-\tau)} \norm{v(\tau)\otimes v(\tau)}_{\frac{d}{1+\delta}} d\tau\cr
&&+ \int_0^t [h_d(t-\tau)]^{\frac{1+\delta}{d}}e^{-\left(d-1  + \frac{\gamma_{d,d} + \gamma_{d/(1+\delta),d}}{2} \right)(t-\tau)} \left\| F(\tau) \right\|_{\frac{d}{1+\delta}} d\tau\cr
&\le& \| u(0)\|_{d} + \int_0^t [h_d(t-\tau)]^{\frac{1+\delta}{d}}[h_d(\tau)]^{\frac{1-\delta}{d}}e^{-\hat{\beta}_1(t-\tau)}\norm{v(\tau)}_{\frac{d}{\delta}}\norm{v(\tau)}_d d\tau\cr
&&+\int_0^t C^{ \frac{1+\delta}{d}} \left((t-\tau)^{-\frac{1+\delta}{2}}+1\right) e^{-\hat{\beta}_1 (t-\tau)} d\tau \|F\|_{\infty,\frac{d}{1+\delta}}\cr
&\le& \| u(0)\|_{d} + \norm{v}^2_{\mathcal{X}}\int_0^t [h_d(t-\tau)]^{\frac{1+\delta}{d}}[h_d(\tau)]^{\frac{1-\delta}{d}}e^{-\hat{\beta}_1(t-\tau)}e^{-\sigma\tau}d\tau\cr
&&+\int_0^t C^{ \frac{1+\delta}{d}} \left((t-\tau)^{-\frac{1+\delta}{2}}+1\right) e^{-\hat{\beta}_1 (t-\tau)} d\tau \|F\|_{\infty,\frac{d}{1+\delta}}\cr
&\le& \|u(0)\|_{d} + M_2\norm{v}^2_{\mathcal{X}} + \tilde{C}\left(\hat{\beta}^{\theta-1}_1{\mathbf{\Gamma}}(1-\theta) + \frac{1}{\hat{\beta}_1} \right)P,
\end{eqnarray}
where $\hat{\beta}_1 = d-1  + \frac{\gamma_{d,d} + \gamma_{d/(1+\delta),d}}{2}$, $\tilde{C}=C^{\frac{\delta}{d}+\frac{1}{d}}>0$, $0<\theta=\dfrac{1+\delta}{2}<1$, ${\mathbf{\Gamma}}$ is the gamma function and, the integral 
$$\int_0^t [h_d(t-\tau)]^{\frac{1+\delta}{d}}[h_d(\tau)]^{\frac{1-\delta}{d}}e^{-\hat{\beta}_1(t-\tau)}e^{-\sigma\tau}d\tau \leqslant M_2 <+\infty$$
converges by the mean of gamma and beta functions (see Appendix).

Finally, we estimate $[h_d(t)]^{-\frac{1-\delta}{d}}e^{\sigma t}\norm{u(t)}_{\frac{d}{\delta}}$.
By using Assertion iii) in Lemma \ref{estimates}, we have
\begin{eqnarray}\label{part31}
&&[h_d(t)]^{-\frac{1-\delta}{d}}e^{\sigma t}\norm{\int_0^t e^{-(t-\tau)\A} \p\dive (-v(\tau)\otimes v(\tau)) d\tau}_{\frac{d}{\delta}} \cr
&\leqslant& C^{\frac{1-\delta}{d}}e^{\sigma t}\int_0^t \norm{e^{-(t-\tau)\A} \p\dive (-v(\tau)\otimes v(\tau))}_{\frac{d}{\delta}} d\tau \hbox{   (because   } [h_d(t)]^{-\frac{1-\delta}{d}}\leqslant C^{\frac{1-\delta}{d}})\cr
&\leqslant& C^{\frac{1-\delta}{d}}e^{\sigma t}\int_0^t [h_d(t-\tau)]^{\frac{1+\delta}{d}} e^{-{\beta}(t-\tau)} \|v(\tau)\otimes v(\tau)\|_{\frac{d}{2\delta}}d\tau\cr
&\leqslant& C^{\frac{1-\delta}{d}}e^{\sigma t}\int_0^t [h_d(t-\tau)]^{\frac{1+\delta}{d}} e^{-{\beta}(t-\tau)} \|v(\tau)\|_{\frac{d}{\delta}} \|v(\tau)\|_{\frac{d}{\delta}}d\tau\cr
&\leqslant& \norm{v}^2_{\mathcal{X}} C^{\frac{1-\delta}{d}}\int_0^t [h_d(t-\tau)]^{\frac{1+\delta}{d}}[h_d(\tau)]^{\frac{2(1-\delta)}{d}}e^{-({\beta}-\sigma)(t-\tau)} e^{-\sigma\tau}d\tau\cr
&\leqslant& M_3\norm{v}^2_{\mathcal{X}},
\end{eqnarray}
where ${\beta} = d-1 + \frac{\gamma_{d/\delta,d/\delta} + \gamma_{d/(2\delta),d/\delta}}{2}$ and
$$C^{\frac{1-\delta}{d}}\int_0^t [h_d(t-\tau)]^{\frac{1+\delta}{d}}[h_d(\tau)]^{\frac{2(1-\delta)}{d}}e^{-({\beta}-\sigma)(t-\tau)}e^{-\sigma\tau}d\tau \leq M_3 <+\infty$$
by the mean of gamma and beta functions (see Appendix). On the other hand, we have
\begin{eqnarray}\label{part32}
&&[h_d(t)]^{-\frac{1-\delta}{d}}e^{\sigma t}\norm{\int_0^t e^{-(t-\tau)\A} \p\dive F(\tau) d\tau}_{\frac{d}{\delta}} \cr
&\leqslant& C^{\frac{1-\delta}{d}}e^{\sigma t}\int_0^t \norm{e^{-(t-\tau)\A} \p\dive F(\tau)}_{\frac{d}{\delta}}d\tau \hbox{   (because   } [h_d(t)]^{-\frac{1-\delta}{d}}\leqslant C^{\frac{1-\delta}{d}})\cr
&\leqslant& C^{\frac{1-\delta}{d}}e^{\sigma t}\int_0^t [h_d(t-\tau)]^{\frac{1+\delta}{d}} e^{-{\beta}(t-\tau)} \|F(\tau)\|_{\frac{d}{2\delta}}d\tau\cr
&\leqslant& \sup_{t>0}[h_d(t)]^{\frac{1-\delta}{d}}e^{\sigma t}\|F(t)\|_{\frac{d}{2\delta}}C^{\frac{1-\delta}{d}}\int_0^t [h_d(t-\tau)]^{\frac{1+\delta}{d}}[h_d(\tau)]^{\frac{1-\delta}{d}}e^{-({\beta}-\sigma)(t-\tau)} d\tau \cr
&\leqslant& M_4P,
\end{eqnarray}
where
$$C^{\frac{1-\delta}{d}}\int_0^t [h_d(t-\tau)]^{\frac{1+\delta}{d}}[h_d(\tau)]^{\frac{1-\delta}{d}}e^{-({\beta}-\sigma)(t-\tau)} d\tau \leqslant M_4 <+\infty.$$
Combining inequalities \eqref{part31} and \eqref{part32} with noting that $[h_d(t)]^{-\frac{1-\delta}{d}}e^{\sigma t}\norm{u(0)}_{\frac{d}{\delta}} \leqslant \norm{u(0)}_d $, we get
\begin{equation}\label{Part3}
[h_d(t)]^{-\frac{1-\delta}{d}}e^{\sigma t} \norm{u(t)}_{\frac{d}{\delta}} \leqslant \| u(0)\|_d + M_3\norm{v}^2_{\mathcal{X}}+ M_4P.
\end{equation}

From inequalities \eqref{Part1}, \eqref{Part2} and \eqref{Part3}, we obtain
\begin{equation*}
\norm{u(t)}_{\mathcal{X}} \leq \norm{u(0)}_{L^p\cap L^d} + M \norm{v}^2_{\mathcal{X}} + NP,
\end{equation*}
where $M=M_1+M_2+M_3$ and $N=\tilde{C}\left(\beta^{\theta-1}_1{\mathbf{\Gamma}}(1-\theta) + \frac{1}{\beta_1} \right) + \tilde{C}\left(\hat{\beta}_1^{\theta-1}{\mathbf{\Gamma}}(1-\theta) + \frac{1}{\hat{\beta}_1} \right)+ M_4$. This is Inequality \eqref{CoreEstimates} and our proof is completed.

\noindent
\item[ii)] Similar to the proof of Assertion $i)$, we need only to prove the boundedness of $\| u\|_{\mathbb{X}}$. Also denote that $\norm{\cdot}_{C_b(\r,L^r(Z))} = \norm{\cdot}_{\infty,r}$ for $r>0$ and $Z$ is $\Gamma(T\mathcal{M})$ or $\Gamma(T\mathcal{M}\otimes T\mathcal{M})$.
This boundedness is proven by the similar estimates as in the proof of Assertion i), with noting that we have not the initial data norm $\norm{u(0)}_{L^p\cap L^d}$ and the integrals in the right hand-sides of inequalities \eqref{Part1}, \eqref{Part2}, \eqref{part31} and \eqref{part32} are taken on $(-\infty,\, t]$ (see the boundedness of these integrals in Appendix). In particular,
\begin{eqnarray}\label{PPPart1}
\|u(t)\|_{p} &\le& \int_{-\infty}^t \|e^{-(t-\tau)\A}\dive (-v(\tau)\otimes v(\tau)+ F(\tau))\|_{p}d\tau\cr 
&\le& \int_{-\infty}^t [h_d(t-\tau)]^{\frac{\delta}{d}+\frac{1}{d}}e^{-\left(d-1  + \frac{\gamma_{p,p} + \gamma_{dp/(1+\delta p),p}}{2} \right)(t-\tau)} \norm{v(\tau)\otimes v(\tau)}_{\frac{dp}{d+\delta p}} d\tau\cr
&&+ \int_{-\infty}^t [h_d(t-\tau)]^{\frac{\delta}{d}+\frac{1}{d}}e^{-\left(d-1  + \frac{\gamma_{p,p} + \gamma_{dp/(1+\delta p),p}}{2} \right)(t-\tau)} \left\| F(\tau) \right\|_{\frac{dp}{d+\delta p}} d\tau\cr
&\le& \int_{-\infty}^t [h_d(t-\tau)]^{\frac{1+\delta}{d}}e^{-\beta_1(t-\tau)}\norm{v(\tau)}_{\frac{d}{\delta}}\norm{v(\tau)}_p d\tau\cr
&&+\int_{-\infty}^t C^{ \frac{\delta}{d}+\frac{1}{d}} \left((t-\tau)^{-\frac{1+\delta}{2}}+1\right) e^{-\beta_1 (t-\tau)} d\tau \|F\|_{\infty,\frac{dp}{d+\delta p}}\cr
&\le& \norm{v}^2_{\mathbb{X}}\int_{-\infty}^t [h_d(t-\tau)]^{\frac{1+\delta}{d}}\lambda^{-1}(\tau)e^{-\beta_1(t-\tau)}d\tau\cr
&&+\int_{-\infty}^t C^{ \frac{1+\delta}{d}} \left((t-\tau)^{-\frac{1+\delta}{2}}+1\right) e^{-\beta_1 (t-\tau)} d\tau \|F\|_{\infty,\frac{dp}{d+\delta p}}\cr
&\le& \tilde{M}_1\norm{v}^2_{\mathbb{X}} + \tilde{C}\left(\beta_1^{\theta-1}{\mathbf{\Gamma}}(1-\theta) + \frac{1}{\beta_1} \right)\tilde{P},
\end{eqnarray}
where  
$$\int_{-\infty}^t [h_d(t-\tau)]^{\frac{1+\delta}{d}}\lambda^{-1}(\tau)e^{-\beta_1(t-\tau)}d\tau \leqslant \tilde{M}_1 <+\infty.$$
Similarly, we can estimate that
\begin{eqnarray}\label{PPPart2}
\|u(t)\|_{d} &\le& \int_{-\infty}^t \|e^{-(t-\tau)\A}\dive (-v(\tau)\otimes v(\tau)+ F(\tau))\|_{d}d\tau\cr 
&\le& \norm{v}^2_{\mathbb{X}}\int_{-\infty}^t [h_d(t-\tau)]^{\frac{1+\delta}{d}}\lambda^{-1}(\tau)e^{-\hat{\beta}_1(t-\tau)}d\tau\cr
&&+\int_{-\infty}^t C^{ \frac{1+\delta}{d}} \left((t-\tau)^{-\frac{1+\delta}{2}}+1\right) e^{-\hat{\beta}_1 (t-\tau)} d\tau \|F\|_{\infty,\frac{d}{1+\delta}}\cr
&\le& \tilde{M}_2\norm{v}^2_{\mathbb{X}} + \tilde{C}\left(\hat{\beta}_1^{\theta-1}{\mathbf{\Gamma}}(1-\theta) + \frac{1}{\hat{\beta}_1} \right)\tilde{P},
\end{eqnarray}
where 
$$\int_{-\infty}^t [h_d(t-\tau)]^{\frac{1+\delta}{d}}\lambda^{-1}(\tau)e^{-\hat{\beta}_1(t-\tau)}d\tau \leqslant \tilde{M}_2 <+\infty.$$
In order to estimate $\lambda(t)\norm{u(t)}_{\frac{d}{\delta}}$, we have
\begin{eqnarray}\label{pppart31}
&&\lambda(t)\norm{\int_{-\infty}^t e^{-(t-\tau)\A} \p\dive (-v(\tau)\otimes v(\tau)) d\tau}_{\frac{d}{\delta}} \cr
&\leqslant& \lambda(t)\int_{-\infty}^t \norm{e^{-(t-\tau)\A} \p\dive (-v(\tau)\otimes v(\tau))}_{\frac{d}{\delta}} d\tau \cr
&\leqslant& \lambda(t)\int_{-\infty}^t [h_d(t-\tau)]^{\frac{1+\delta}{d}} e^{-{\beta}(t-\tau)} \|v(\tau)\otimes v(\tau)\|_{\frac{d}{2\delta}}d\tau\cr
&\leqslant& \lambda(t)\int_{-\infty}^t [h_d(t-\tau)]^{\frac{1+\delta}{d}} e^{-{\beta}(t-\tau)} \|v(\tau)\|_{\frac{d}{\delta}} \|v(\tau)\|_{\frac{d}{\delta}}d\tau\cr
&\leqslant& \norm{v}^2_{\mathbb{X}} \left(\lambda(t)\int_{-\infty}^t [h_d(t-\tau)]^{\frac{1+\delta}{d}}\lambda^{-2}(\tau)e^{-{\beta}(t-\tau)} d\tau \right)\cr
&\leqslant& \tilde{M}_3\norm{v}^2_{\mathbb{X}}.
\end{eqnarray}
where 
$$\lambda(t)\int_{-\infty}^t [h_d(t-\tau)]^{\frac{1+\delta}{d}}\lambda^{-2}(\tau)e^{-{\beta}(t-\tau)}d\tau \leq \tilde{M}_3 <+\infty.$$
Moreover,
\begin{eqnarray}\label{pppart32}
&&\lambda(t)\norm{\int_{-\infty}^t e^{-(t-\tau)\A} \p\dive F(\tau) d\tau}_{\frac{d}{\delta}} \cr
&\leqslant& \lambda(t)\int_{-\infty}^t \norm{e^{-(t-\tau)\A} \p\dive F(\tau)}_{\frac{d}{\delta}}d\tau \cr
&\leqslant& \lambda(t)\int_{-\infty}^t [h_d(t-\tau)]^{\frac{1+\delta}{d}} e^{-{\beta}(t-\tau)} \|F(\tau)\|_{\frac{d}{2\delta}}d\tau\cr
&\leqslant& \sup_{t\in \mathbb{R}}\lambda(t)\|F(t)\|_{\frac{d}{2\delta}} \left(\lambda(t)\int_{-\infty}^t [h_d(t-\tau)]^{\frac{1+\delta}{d}}\lambda^{-1}(\tau) e^{-{\beta}(t-\tau)} d\tau \right)\cr
&\leqslant& \tilde{M}_4\tilde{P},
\end{eqnarray}
where
$$\lambda(t)\int_{-\infty}^t [h_d(t-\tau)]^{\frac{1+\delta}{d}}\lambda^{-1}(\tau)e^{-{\beta}(t-\tau)} d\tau \leqslant \tilde{M}_4 <+\infty.$$
Inequalities \eqref{pppart31} and \eqref{pppart32} lead to
\begin{equation}\label{PPPart3}
\lambda(t)\norm{u(t)}_{\frac{d}{\delta}} \leqslant \tilde{M}_3\norm{v}^2_{\mathbb{X}}+ \tilde{M}_4\tilde{P}.
\end{equation}

Since inequalities \eqref{PPPart1}, \eqref{PPPart2} and \eqref{PPPart3}, we obtain that
\begin{equation*}
\norm{u(t)}_{\mathbb{X}} \leq \tilde{M} \norm{v}^2_{\mathbb{X}} + \tilde{N}\tilde{P},
\end{equation*}
where $\tilde{M} = \tilde{M}_1+\tilde{M}_2+\tilde{M}_3$ and $\tilde{N} = \tilde{C}\left(\beta_1^{\theta-1}{\mathbf{\Gamma}}(1-\theta) + \frac{1}{\beta_1} \right) + \tilde{C}\left(\hat{\beta}_1^{\theta-1}{\mathbf{\Gamma}}(1-\theta) + \frac{1}{\hat{\beta}_1} \right)+ \tilde{M}_4$.
Therefore, the boundedness of $u$ holds.
  
\end{proof} 
\begin{remark}
In Lemma \ref{Thm:linear}, the boundedness of mild solution $u(t)$ defined on half time-line axis (resp. full time-line axis) is proved in the Banach space $\mathcal{X}= C_b(\mathbb{R}_+,(L^p\cap L^d\cap L^{\frac{d}{\delta}})(\Gamma(T\mathcal{M})))$ (resp. $\mathbb{X}=C_b(\mathbb{R},(L^p\cap L^d\cap L^{\frac{d}{\delta}})(\Gamma(T\mathcal{M})))$) with the mixed norm \eqref{space1} (resp. \eqref{space2}). This complication comes from the fact that $1<p\leq d$. By this condition we can not prove the boundedness of $u(t)$ in the single space $L^p(\Gamma(T\mathcal{M}))$ by the mean of gamma functions as well as in our recent work \cite{XVQ2023} for the case $p>d$. 
\end{remark}

For the sake of simplicity, we denote $X= (L^{\frac{dp}{d+\delta p}}\cap L^{\frac{d}{1+\delta}}\cap L^{\frac{d}{2\delta}})(\Gamma(T\mathcal{M}\otimes T\mathcal{M}))$ and $Y= (L^p\cap L^d \cap L^{\frac{d}{\delta}})(\Gamma(T\mathcal{M}))$. By using Lemma \ref{Thm:linear}, we can define the solution operator $S: C_b(\r_+, X\times Y)\to \mathcal{X}=C_b(\r_+,Y)$ of Equation \eqref{CauchyHeat} as follows
\begin{align*}
S: C_b(\r_+,X\times Y) &\rightarrow \mathcal{X}=C_b(\r_+,Y)\cr
(F,v)&\mapsto S(F),
\end{align*}
where $X\times Y$ is the Cartesian product space equipped with the norm $\norm{\cdot}_{X\times Y} = \norm{\cdot}_X + \norm{\cdot}_Y$ and
\begin{equation}\label{SolOpe}
S(F,v)(t):= u(0)+\int_0^t e^{-(t-\tau)\A}\dive[-v(\tau)\otimes v(\tau)+F(\tau)]d\tau.
\end{equation}
Here, we define the norm of Banach space $C_b(\mathbb{R}_+,X\times Y)$ by
$$\norm{(F,v)}_{C_b(\r_+,X\times Y)} = \sup_{t>0}\left(\norm{F(t)}_{\frac{dp}{d+\delta p}} + \norm{F(t)}_{\frac{d}{1+\delta}} + [h_d(t)]^{-\frac{1-\delta}{d}}e^{\sigma t}\norm{F(t)}_{\frac{d}{2\delta}} \right) + \norm{v}_{\mathcal{X}}.$$
The norm of space of AAP- functions $AAP(\r_+,X\times Y)$ is inherited from the above norm.

In order to prove the below theorem, we also define the norm of Banach space $C_b(\mathbb{R},X\times Y)$ by
$$\norm{(F,v)}_{C_b(\r,X\times Y)} = \sup_{t\in \mathbb{R}}\left(\norm{F(t)}_{\frac{dp}{d+\delta p}} + \norm{F(t)}_{\frac{d}{1+\delta}} + \lambda(t)\norm{F(t)}_{\frac{d}{2\delta}} \right) + \norm{v}_{\mathbb{X}},$$
where $\lambda(t)$ is given by \eqref{lambda}. This norm is used for the space of almost periodic functions $AP(\r,X\times Y)$.

Using Lemma \ref{Thm:linear}, we prove a Massera-type principle to obtain the main theorem of this section. In particular, in Theorem \ref{pest} below, we extend Theorem 3.6 in \cite{XVQ2023} to prove the existence of forward asymptotically almost periodic (AAP-) mild solution of Equation \eqref{CauchyHeat} in the phase space $\mathcal{X}=C_b(\mathbb{R}_+,Y=(L^p\cap L^d\cap L^{\frac{d}{\delta}})(\Gamma(T\mathcal{M})))\, \,(1<p\leq d)$ with the norm given by \eqref{space1}.
\begin{theorem}\label{pest}
Let $(\mathcal{M},g)$ be a $d$-dimensional real hyperbolic manifold with $d \geqslant 2$.
Let $0<\delta<\min\left\{1,\, \dfrac{d}{8} \right\}$, suppose that functions $F$ and $v$ are given such that the function $t \mapsto (F(t),v(t))$ belongs  to $AAP(\mathbb{R}_+,X\times Y)$, where $X= (L^{\frac{dp}{d+\delta p}}\cap L^{\frac{d}{1+\delta}}\cap L^{\frac{d}{2\delta}})(\Gamma(T\mathcal{M}\otimes T\mathcal{M}))$ and $Y= (L^p\cap L^d \cap L^{\frac{d}{\delta}})(\Gamma(T\mathcal{M}))\, (1<p\leq d)$.
Then, Problem \eqref{CauchyHeat} has one and only one forward AAP- mild solution $\hat{u} \in AAP(\mathbb{R}_+,Y)$ satisfying 
\begin{equation}\label{esper}
\|\hat{u}\|_{\mathcal{X}}\le \|u(0)\|_{L^p\cap L^d} + M\norm{v}^2_{\mathcal{X}} + NP,
\end{equation}
where $M,\, N$ and $P$ are given in Assertion i) in Lemma \ref{Thm:linear}. 
\end{theorem} 
\def\xcal{\mathcal X}
\begin{proof} 
By using Assetion i) in Lemma \ref{Thm:linear}, it is enough to show that, the solution operator $S$ maps $AAP(\mathbb{R}_+,X\times Y)$ into $AAP(\mathbb{R}_+,Y)$. 

Indeed, for each $(F,v)\in AAP( \mathbb{R}_+, X\times Y) $, there exist  $(H,\eta) \in AP(\mathbb{R},X\times Y)$ and $ (\Phi,\omega) \in C_0(\r_+,X\times Y)$ such that  $F(t) =  H(t)+ \Phi(t)$ and $v(t)=\eta(t)+\omega(t)$ for all $t \in \r_+$. Using \eqref{SolOpe}, we have
\begin{eqnarray}\label{seperate}
&&S(F,v)(t) = e^{-t\cal{A}}u(0) + \int_{0}^t  e^{-(t-\tau) \mathcal{A}}\mathbb{P}\dive[-v(\tau)\otimes v(\tau)+ F(\tau)] d\tau\cr
&=& e^{-t\cal{A}}u(0) + \int_{0}^t  e^{-(t-\tau) \mathcal{A}}\mathbb{P}\dive H(\tau)d\tau+\int_{0}^t  e^{-(t-\tau) \mathcal{A}}\mathbb{P}\dive \Phi(\tau)d\tau\cr
&&- \int_{0}^t  e^{-(t-\tau) \mathcal{A}}\mathbb{P}\dive [\eta(\tau)\otimes \eta(\tau)]d\tau - \int_{0}^t  e^{-(t-\tau) \mathcal{A}}\mathbb{P}\dive [\omega(\tau)\otimes v(\tau)] d\tau \cr
&&- \int_{0}^t  e^{-(t-\tau) \mathcal{A}}\mathbb{P}\dive [\eta(\tau)\otimes \omega(\tau)] d\tau\cr
&=& e^{-t\cal{A}}u(0) + \int_{-\infty}^t  e^{-(t-\tau) \mathcal{A}}\mathbb{P}\dive H(\tau)d\tau + \int_{0}^t  e^{-(t-\tau) \mathcal{A}}\mathbb{P}\dive \Phi(\tau)d\tau \cr
&&- \int_{-\infty}^0 e^{-(t-\tau)\mathcal{A}}\mathbb{P}\dive H(\tau)d\tau - \int_{-\infty}^t  e^{-(t-\tau) \mathcal{A}}\mathbb{P}\dive[ \eta(\tau)\otimes \eta(\tau)]d\tau \cr
&&-\int_{0}^t  e^{-(t-\tau) \mathcal{A}}\mathbb{P}\dive[\omega(\tau)\otimes v(\tau)]d\tau - \int_{0}^t  e^{-(t-\tau) \mathcal{A}}\mathbb{P}\dive [\eta(\tau)\otimes \omega(\tau)] d\tau\cr
&& + \int_{-\infty}^0 e^{-(t-\tau)\mathcal{A}}\mathbb{P}\dive[\eta(\tau)\otimes \eta(\tau)]d\tau
\end{eqnarray}
for all $t\in \mathbb{R}_+$.

Setting
\begin{eqnarray*}
&&\hat{S}(H)(t)=\int_{-\infty}^te^{-(t-\tau)\mathcal{A}} \mathbb{P} \dive H (\tau) d\tau,\cr
&& \hat{S}(\eta\otimes \eta)(t)=\int_{-\infty}^te^{-(t-\tau)\mathcal{A}} \mathbb{P} \dive [\eta(\tau)\otimes \eta(\tau)]d\tau,
\end{eqnarray*}
and
\begin{eqnarray*}
&&\tilde{S}(\Phi)(t) = \int_{0}^t  e^{-(t-\tau) \mathcal{A}}\mathbb{P}\dive \Phi(\tau)d\tau,\cr
&&\tilde{S}(\omega\otimes v)(t) = \int_{0}^t  e^{-(t-\tau) \mathcal{A}}\mathbb{P}\dive [\omega(\tau)\otimes v(\tau)]d\tau,\cr
&&\tilde{S}(\eta\otimes \omega)(t) = \int_{0}^t  e^{-(t-\tau) \mathcal{A}}\mathbb{P}\dive [\eta(\tau)\otimes \omega(\tau)]d\tau.
\end{eqnarray*}
Following the proof of Lemma \ref{Thm:linear}, the functions $\hat{S}(H)$ and $\hat{S}(\eta\otimes \eta)$ are bounded in $\mathbb{X}=C_b(\mathbb{R},Y)$ with norm given by \eqref{space2}; the functions  $\tilde{S}(\Phi),\,\tilde{S}(\omega\otimes v)$ and $\tilde{S}(\eta\otimes \omega)$ are bounded in $\mathcal{X}=C_b(\mathbb{R}_+,Y)$ with norm given by \eqref{space1}. Since \eqref{seperate}, we have
\begin{eqnarray}
S(F,v)(t) &=& e^{-t\cal{A}}u(0) + \left( \hat{S}(H)(t) - \hat{S}(\eta\otimes \eta)(t) \right)\cr
&&+ \left( \tilde{S}(\Phi)(t) - \tilde{S}(\omega\otimes v)(t) - \tilde{S}(\eta\otimes \omega)(t) \right) \cr
&&- \left( e^{-t\mathcal{A}}\hat{S}(H)(0) - e^{-t\mathcal{A}}\hat{S}(\eta\otimes \eta)(0) \right) \hbox{   for all  } t\in\mathbb{R_+}.
\end{eqnarray}

Now, we prove that $S(F,v)\in AAP(\r_+,Y)$ by three following steps:

\underline{\bf Step 1:} The function $ \hat{S}(H)-\hat{S}(\eta\otimes \eta)$ belongs to $AP(\mathbb{R}, Y)$. Indeed, since $ (H,\eta)\in AP(\mathbb{R}, X\times Y)$, for each $ \epsilon  > 0$, there exists $l_{\epsilon}>0 $ such that every interval of length $l_{\epsilon}$ contains at least a number $T $ with the following property
$$\sup_{t \in \r } \left( \| H(t+T)  - H(t) \|_{C_b(\r,X)} + \| \eta(t+T) - \eta(t) \|_{C_b(\r,Y)} \right) < \epsilon,$$
where
$$\norm{H}_{C_b(\r,X)}= \sup_{t\in \mathbb{R}} \left(\norm{H(t)}_{\frac{dp}{d+\delta p}} + \norm{H(t)}_{\frac{d}{1+\delta}} + \lambda(t)\norm{H(t)}_{\frac{d}{2\delta}} \right)$$
and $\norm{\eta}_{C_b(\r,Y)}= \norm{\eta}_{\mathbb{X}}$ given by \eqref{space2}.

By the same way as in the proof of Assertion ii) of Lemma \ref{Thm:linear}, we can estimate
\begin{eqnarray}\label{AP1}
&&\left\|\hat{S}(H)(t+T) - \hat{S}(H)(t)\right\|_Y \cr
&=& \left\|\int_{-\infty}^{t+T} e^{-(t+T-\tau)\mathcal{A}}\mathbb{P}[\dive H(\tau)] d\tau - \int_{-\infty}^{t} e^{-(t-\tau)\mathcal{A}}\mathbb{P}[\dive H(\tau)] d\tau \right\|_Y \cr
&=& \left\|\int_{0}^\infty e^{-\tau\mathcal{A}}\mathbb{P}[\dive H(t+T-\tau)-\dive H(t-\tau)] d\tau \right\|_Y \cr
&\leqslant& \tilde{N}\left\| H(\cdot+T) - H(\cdot)\right\|_{C_b(\r,X)} < \tilde{N}\epsilon,
\end{eqnarray}
for all $t \in \r$, where $\tilde{N}$ is determined as in Assertion ii) of Lemma \ref{Thm:linear}. 
Moreover, we also have
\begin{eqnarray}\label{AP2}
&&\left\|\hat{S}(\eta\otimes \eta)(t+T) - \hat{S}(\eta\otimes \eta)(t)\right\|_Y\cr
&=& \left\|\int_{-\infty}^{t+T} e^{-(t+T-\tau)\mathcal{A}}\mathbb{P}\dive[ (\eta\otimes \eta)(\tau)] d\tau - \int_{-\infty}^{t} e^{-(t-\tau)\mathcal{A}}\mathbb{P}\dive[ (\eta\otimes \eta)(\tau)] d\tau\right\|_Y \cr
&=& \left\|\int_{0}^\infty e^{-\tau\mathcal{A}}\mathbb{P}\dive[ (\eta\otimes \eta)(t+T-\tau)-(\eta\otimes \eta)(t-\tau)] d\tau \right\|_Y \cr
&\leqslant& \left\| \int_{0}^\infty  e^{-\tau\mathcal{A}}\mathbb{P}\dive [(\eta(t+T-\tau)- \eta(t-\tau))\otimes \eta(t+T-\tau)] d\tau \right\|_Y\cr
&&+ \left\|\int_{0}^\infty  e^{-\tau\mathcal{A}}\mathbb{P}\dive [\eta(t-\tau)\otimes (\eta(t+T-\tau)- \eta(t-\tau))]d\tau \right\|_{Y} \cr
&\leqslant& 2\tilde{M} \norm{\eta(\cdot+T)- \eta(\cdot)}_{C_b(\r,Y)}\norm{\eta}_{C_b(\r,Y)} < 2\tilde{M}\norm{\eta}_{C_b(\r,Y)}\epsilon,
\end{eqnarray}
for all $t\in \mathbb{R}$, where $\tilde{M}$ is given as in Assertion ii) of Lemma \ref{Thm:linear}.

Combining inequalities \eqref{AP1} and \eqref{AP2}, we obtain that 
$$ \norm{(\hat{S}(H)-\hat{S}(\eta\otimes \eta))(t+T) - (\hat{S}(H)-\hat{S}(\eta\otimes \eta))(t)}_{Y} <(\tilde{N} + 2\tilde{M}\norm{\eta}_{\mathbb{X}})\epsilon.$$
This shows that $\hat{S}(H)-\hat{S}(\eta\otimes \eta)$ belongs to $AP(\mathbb{R}, Y)$.

\underline{\bf Step 2:} We show that $\tilde{S}(\Phi) - \tilde{S}(\omega\otimes v) - \tilde{S}(\eta\otimes \omega)$ belongs to $C_0(\r_+,Y)$. Indeed, the first term $\tilde{S}(\Phi)(t)$ can be rewritten as follows
\begin{eqnarray*}
\tilde{S}(\Phi)(t) &=& \int_{0}^t  e^{-(t-\tau) \mathcal{A}}\mathbb{P}\dive \Phi(\tau)d\tau\hbox{   } (t\in\mathbb{R_+}) \cr
&=&\int_{0}^{t/2}  e^{-(t-\tau) \mathcal{A}}\mathbb{P}\dive \Phi(\tau)d\tau +\int_{t/2}^t  e^{-(t-\tau) \mathcal{A}}\mathbb{P}\dive \Phi(\tau)d\tau \cr
&=& S_1(\Phi)(t) +S_2(\Phi)(t).
\end{eqnarray*}
By the same way to establish inequalities \eqref{Part1}, \eqref{Part2}, \eqref{part31} and \eqref{part32} in the proof of Assertion i) in Lemma \ref{Thm:linear}, for $t>2$, we obtain
\begin{eqnarray*}
&&\norm{S_1(\Phi)(t)}_p \leqslant \int_0^{t/2} C^{\frac{1+\delta}{d}} \left( (t-\tau)^{-\frac{1+\delta}{2}}+1\right)e^{-\beta_1(t-\tau)}d\tau \norm{\Phi}_{\infty,\frac{dp}{d+\delta p}}\cr
&\leqslant& C^{\frac{1+\delta}{d}} \left( \left( \frac{t}{2}\right)^{-\frac{1+\delta}{2}}+1\right)\frac{1}{\beta_1}(e^{-\frac{\beta_1 t}{2}} - e^{-\beta_1 t}) \norm{\Phi}_{C_b(\r_+,X)} \to 0, \hbox{   when  } t \to +\infty,\cr
&&\norm{S_1(\Phi)(t)}_d \leqslant \int_0^{t/2} C^{\frac{1+\delta}{d}} \left( (t-\tau)^{-\frac{1+\delta}{2}}+1\right)e^{-\hat{\beta}_1(t-\tau)}d\tau \norm{\Phi}_{\infty,\frac{d}{1+\delta}}\cr
&\leqslant& C^{\frac{1+\delta}{d}} \left( \left( \frac{t}{2}\right)^{-\frac{1+\delta}{2}}+1\right)\frac{1}{\hat\beta}_1(e^{-\frac{\hat{\beta}_1 t}{2}} - e^{-\hat{\beta}_1 t}) \norm{\Phi}_{C_b(\r_+,X)} \to 0, \hbox{   when  } t \to +\infty,\cr
&&[h_d(t)]^{-\frac{1-\delta}{d}} e^{\sigma t}\norm{S_1(\Phi)(t)}_{\frac{d}{\delta}} \cr
&\leqslant& C^{\frac{1-\delta}{d}}\int_0^{t/2} [h_d(t-\tau)]^{\frac{1+\delta}{d}}[h_d(\tau)]^{\frac{1-\delta}{d}}e^{-({\beta}-\sigma)(t-\tau)}d\tau \sup_{t>0}[h_d(t)]^{\frac{1-\delta}{d}}e^{\sigma t}\norm{\Phi(t)}_{\frac{d}{2\delta}}\cr
&\leqslant& C^{\frac{1-\delta}{d}}\left( \left( \frac{t}{2} \right)^{-\frac{1+\delta}{2}}+1\right)\int_0^{t/2} [h_d(\tau)]^{\frac{1-\delta}{d}}e^{-({\beta}-\sigma)(t-\tau)}d\tau \norm{\Phi}_{C_b(\r_+,X)}\cr
&\leqslant& C^{\frac{1-\delta}{d}}\left( \left( \frac{t}{2} \right)^{-\frac{1+\delta}{2}}+1\right)\left(\int_0^1 \tau^{-\frac{1-\delta}{2}}e^{-({\beta}-\sigma)(t-\tau)}d\tau + \int_1^{t/2} e^{-({\beta}-\sigma)(t-\tau)}d\tau \right)\norm{\Phi}_{C_b(\r_+,X)} \cr
&\leqslant& C^{\frac{1-\delta}{d}}\left( \left( \frac{t}{2} \right)^{-\frac{1+\delta}{2}}+1\right)\left( e^{-({\beta}-\sigma)(t-1)} - \frac{1}{{\beta}-\sigma}(e^{-({\beta}-\sigma)(t-1)}-e^{-({\beta}-\sigma)\frac{t}{2}}) \right) \norm{\Phi}_{C_b(\r_+,X)}\cr
&& \to 0, \hbox{   when  } t \to +\infty.
\end{eqnarray*}
These limits lead to
\begin{equation}\label{lim11}
\lim_{t\to +\infty}\norm{S_1(\Phi)(t)}_{Y}=0.
\end{equation}
On the other hand, we have $\lim\limits_{t\to +\infty}\norm{\Phi(t)}_X = 0$, then for all $\epsilon>0$, there exists $t_0>0$ large enough such that for all $t>t_0$, we have $\norm{\Phi(t)}_X<\epsilon$.
Using this and again by the same way to prove inequalities \eqref{Part1}, \eqref{Part2}, \eqref{part31} and \eqref{part32}, we can show that
\begin{eqnarray*}
\norm{S_2(\Phi)(t)}_{Y} \leqslant N\epsilon \hbox{   for  } t>t_0,
\end{eqnarray*}
where the constant $N$ is given in the proof of Assertion i) of Lemma \ref{Thm:linear}. This yeilds 
\begin{equation}\label{lim12}
\lim_{t\to+\infty}\norm{S_2(\Phi)(t)}_{Y}=0.
\end{equation}
Combining limits \eqref{lim11} and \eqref{lim12}, we get
\begin{equation}\label{lim1}
\lim_{t\to+\infty}\norm{\tilde{S}(\Phi)(t)}_{Y}=0,
\end{equation}
which means that $\tilde{S}(\Phi)\in C_0(\r_+,Y)$.

The proofs of $\tilde{S}(\omega\otimes v)\in C_0(\r_+,Y)$ and $\tilde{S}(\eta\otimes \omega)\in C_0(\r_+,Y)$ are similar, we prove only that $\tilde{S}(\omega\otimes v)\in C_0(\r_+,Y)$. Indeed, we have 
\begin{eqnarray*}
\tilde{S}(\omega\otimes v)(t) &=& \int_{0}^t  e^{-(t-\tau) \mathcal{A}}\mathbb{P}\dive [\omega(\tau)\otimes v(\tau)]d\tau\hbox{   } (t\in\mathbb{R_+}) \cr
&=&\int_{0}^{t/2}  e^{-(t-\tau) \mathcal{A}}\mathbb{P}\dive [\omega(\tau)\otimes v(\tau)]d\tau +\int_{t/2}^t  e^{-(t-\tau) \mathcal{A}}\mathbb{P}\dive [\omega(\tau)\otimes v(\tau)]d\tau \cr
&=& S_3(\omega\otimes v)(t) +S_4(\omega\otimes v)(t).
\end{eqnarray*}
By the same manner to establish inequalities \eqref{Part1}, \eqref{Part2}, \eqref{part31} and \eqref{part32} in the proof of Assertion i) in Lemma \ref{Thm:linear}, for $t>2$, we obtain
\begin{eqnarray*}
&&\|S_3(\omega\otimes v)(t)\|_{p} \leqslant \int_0^{t/2} \|e^{-(t-\tau)\A}\dive [\omega(\tau)\otimes v(\tau)]\|_{p}d\tau\cr 
&\le& \int_0^{t/2} [h_d(t-\tau)]^{\frac{1+\delta}{d}}[h_d(\tau)]^{\frac{1-\delta}{d}}e^{-\beta_1(t-\tau)}e^{-\sigma\tau}d\tau \norm{v}_{\mathcal{X}}\norm{\omega}_{\mathcal{X}}\cr
&\le& \int_0^1 [h_d(\tau)]^{\frac{1-\delta}{d}}e^{-\beta_1(t-\tau)}e^{-\sigma\tau}d\tau \norm{v}_{\mathcal{X}}\norm{\omega}_{\mathcal{X}} + \int_1^{t/2} e^{-\beta_1(t-\tau)}e^{-\sigma\tau}d\tau \norm{v}_{\mathcal{X}}\norm{\omega}_{\mathcal{X}}\cr
&\le& \int_0^1 \tau^{-\frac{1-\delta}{2}}e^{-\beta_1(t-\tau)}e^{-\sigma\tau}d\tau \norm{v}_{\mathcal{X}}\norm{\omega}_{\mathcal{X}} +\int_1^{t/2} e^{-\beta_1(t-\tau)}e^{-\sigma\tau}d\tau \norm{v}_{\mathcal{X}}\norm{\omega}_{\mathcal{X}}\cr
&\le& \left[ \frac{2}{1+\delta}e^{-\beta_1 (t-1)} + \frac{1}{\beta_1}\left( e^{-\frac{\beta_1 t}{2}} - e^{-\beta_1(t-1)} \right) \right] \norm{v}_{\mathcal{X}}\norm{\omega}_{\mathcal{X}} \to 0, \hbox{   when   } t \to +\infty,\cr
&&\|S_3(\omega\otimes v)(t)\|_{d} \leqslant \int_0^{t/2} \|e^{-(t-\tau)\A}\dive [\omega(\tau)\otimes v(\tau)]\|_{d}d\tau\cr 
&\le& \left[ \frac{2}{1+\delta}e^{-\hat{\beta}_1 (t-1)} +  \frac{1}{\hat{\beta}_1}\left( e^{-\frac{\hat{\beta}_1 t}{2}} - e^{-\hat{\beta}_1(t-1)} \right) \right] \norm{v}_{\mathcal{X}}\norm{\omega}_{\mathcal{X}} \to 0, \hbox{   when   } t \to +\infty,\cr
&&[h_d(t)]^{-\frac{1-\delta}{d}}e^{\sigma t}\norm{S_3(\omega\otimes v)(t)}_{\frac{d}{\delta}} \leqslant C^{\frac{1-\delta}{d}}e^{\sigma t}\int_0^{t/2} \norm{e^{-(t-\tau)\A} \p\dive [\omega(\tau)\otimes v(\tau)]}_{\frac{d}{\delta}} d\tau\cr
&\leqslant& C^{\frac{1-\delta}{d}}\int_0^{t/2} [h_d(t-\tau)]^{\frac{1+\delta}{d}}[h_d(\tau)]^{\frac{2(1-\delta)}{d}}e^{-({\beta}-\sigma)(t-\tau)} e^{-\sigma\tau}d\tau \norm{\omega}_{\mathcal{X}}\norm{v}_{\mathcal{X}}\cr
&\leqslant& C^{\frac{1-\delta}{d}}\norm{\omega}_{\mathcal{X}}\norm{v}_{\mathcal{X}} \left( \int_0^{1} \tau^{-(1-\delta)}e^{-({\beta}-\sigma)(t-\tau)} e^{-\sigma\tau}d\tau +\int_1^{t/2} e^{-({\beta}-\sigma)(t-\tau)} e^{-\sigma\tau}d\tau \right)\cr
&\leqslant& C^{\frac{1-\delta}{d}}\norm{\omega}_{\mathcal{X}}\norm{v}_{\mathcal{X}} \left[ \frac{1}{\delta}e^{-({\beta}-\sigma)(t-1)} + \frac{1}{{\beta}-\sigma}\left( e^{-\frac{({\beta}-\sigma) t}{2}} - e^{-({\beta}-\sigma)(t-1)} \right) \right]\cr
&&\to 0, \hbox{   when   } t \to +\infty.
\end{eqnarray*}
These limits lead to
\begin{equation}\label{lim21}
\lim_{t\to +\infty}\norm{S_3(\omega\otimes v)(t)}_{Y}=0.
\end{equation}
On the other hand, we have $\lim\limits_{t\to +\infty}\norm{\omega(t)}_X = 0$, then for all $\epsilon>0$, there exists $t_0>0$ large enough such that for all $t>t_0$, we have $\norm{\omega(t)}_X<\epsilon$.
Using this and again by the same manner to prove inequalities \eqref{Part1}, \eqref{Part2}, \eqref{part31} and \eqref{part32}, we can show that
\begin{eqnarray*}
\norm{S_4(\omega\otimes v)(t)}_{Y} \leqslant N\epsilon \hbox{   for  } t>t_0,
\end{eqnarray*}
where the constant $N$ is given in the proof of Assertion i) of Lemma \ref{Thm:linear}. This yeilds that
\begin{equation}\label{lim22}
\lim_{t\to+\infty}\norm{S_4(\omega\otimes v)(t)}_{Y}=0.
\end{equation}
Combining limits \eqref{lim21} and \eqref{lim22}, we get
\begin{equation}\label{lim2}
\lim_{t\to+\infty}\norm{\tilde{S}(\omega\otimes v)(t)}_{Y}=0,
\end{equation}
which means $\tilde{S}(\omega\otimes v)\in C_0(\r_+,Y)$. Similarly, we have $\tilde{S}(\eta\otimes\omega)\in C_0(\r_+,Y)$.

\underline{\bf Step 3:} Finally, we prove that the function $t\mapsto e^{-t\mathcal{A}}\left(u(0) - \hat{S}(H)(0) + \hat{S}(\eta\otimes \eta)(0)\right)$ belongs to $C_0(\r_+,Y)$. Thanks to Assertion $i)$ in Lemma \ref{estimates}, we have
\begin{eqnarray}
&&\left\| e^{-t\mathcal{A}}\left(u(0)- \hat{S}(H)(0) + \hat{S}(\eta\otimes \eta)(0)\right) \right\|_p\cr
&\leqslant& e^{- \beta_2 t}\left( \left\|u(0) \right\|_p + \norm{\hat{S}(H)(0)}_p + \norm{\hat{S}(\eta\otimes \eta)(0)}_p\right)\cr
&\leqslant&e^{- \beta_2 t}\left( \left\|u(0) \right\|_p + N\norm{H}_{C_b(\r_+,X)} + M\norm{\eta}^2_{\mathcal{X}}\right)  \to 0 \hbox{   when   } t\to +\infty,\label{Term1}\\
&&\left\| e^{-t\mathcal{A}}\left(u(0)- \hat{S}(H)(0) + \hat{S}(\eta\otimes \eta)(0)\right) \right\|_d\cr
&\leqslant& e^{- \hat{\beta}_2 t} \left( \left\|u(0) \right\|_d + N\norm{H}_{C_b(\r_+,X)} + M\norm{\eta}^2_{\mathcal{X}}\right)  \to 0 \hbox{   when   } t\to +\infty,\label{Term2}\\
&&\left\| e^{-t\mathcal{A}}\left(u(0) - \hat{S}(H)(0) + \hat{S}(\eta\otimes \eta)(0)\right) \right\|_{\frac{d}{\delta}} \cr
&\leqslant& [h_d(t)]^{\frac{1-\delta}{d}} e^{- \tilde{\beta}_2 t} \left\| u(0) - \hat{S}(H)(0) + \hat{S}(\eta\otimes \eta)(0)\right\|_d,\label{TTerm3}
\end{eqnarray}
where $ \beta_2 = d-1  + \gamma_{p,p}>0$, $ \hat{\beta}_2 = d-1  + \gamma_{d,d}>0$ and $\tilde{\beta}_2 = d-1 + \gamma_{d,d/\delta}$. By detailed calculation on the formula of $\gamma_{p,q}$ in Lemma \ref{estimates}, we can show that $0<{\beta}<\tilde{\beta}_2$ if $0<\delta<\dfrac{d}{8}$. Therefore, if $0<\delta<\min \left\{ 1,\, \dfrac{d}{8}\right\}$, then Inequality \eqref{TTerm3} leads to 
\begin{eqnarray}\label{Term3}
&&[h_d(t)]^{-\frac{1-\delta}{d}} e^{\sigma t}\left\| e^{-t\mathcal{A}}\left(u(0) - \hat{S}(H)(0) + \hat{S}(\eta\otimes \eta)(0) \right) \right\|_{\frac{d}{\delta}}\cr
&&\leqslant e^{- ({\beta}-\sigma) t} \left\| u(0) - \hat{S}(H)(0) + \hat{S}(\eta\otimes \eta)(0)\right\|_d\cr
&& \leqslant e^{- ({\beta}-\sigma) t} \left( \norm{u(0)}_d + \norm{\hat{S}(H)(0)}_d + \norm{\hat{S}(\eta\otimes \eta)(0)}_d \right)\cr
&& \leqslant e^{- ({\beta}-\sigma) t} \left( \norm{u(0)}_d + N\norm{H}_{C_b(\r_+,X)} + M\norm{\eta}^2_{\mathcal{X}} \right) \to 0 \hbox{   when   } t\to +\infty,
\end{eqnarray}
where the above constants $M$ and $N$ are given in Assertion i) of Lemma \ref{Thm:linear}. By the limits \eqref{Term1}, \eqref{Term2} and \eqref{Term3}, we obtain that
\begin{equation*}
\lim_{t\to +\infty}\norm{e^{-t\mathcal{A}}\left(u(0) - \hat{S}(H)(0) + \hat{S}(\eta\otimes \eta)(0) \right)}_{Y}=0. 
\end{equation*}

Our proof is completed by combinning three above steps.
\end{proof}

\section{The Navier-Stokes equations: existence, uniqueness and asymptotic behaviour of the small forward AAP- mild solutions}\label{S4}
In this section, we consider the existence and exponentially asymptotic behaviour of forward asymptotically almost periodic (AAP-) mild solutions of the incompressible Navier-Stokes equations \eqref{DivNavierStokes}.

The mild solution of Equation \eqref{DivNavierStokes} satisfies the following integral equation 
\begin{equation}\label{MildS}
u(t) = u(0) + \int_0^t e^{-(t-\tau)\mathcal{A}} 
 \mathbb{P}\dive\left(-u(\tau) \otimes u (\tau) + F (\tau)\right) d\tau,
\end{equation}
where $\mathcal{A}u = -(\overrightarrow{\Delta}u - (d-1) u)$ and  
$\mathbb{P}(v) = v + \mathrm{grad}({\Delta}_g)^{-1} \mathrm{div}(v)$.

Now we state and prove the existence and asymptotic behaviour of AAP- mild solution of Equation \eqref{DivNavierStokes} in the phase space $\mathcal{X}=C_b(\r_+, (L^p\cap L^d\cap L^{\frac{d}{\delta}})(\Gamma(T\mathcal{M})))\,\,\, (1<p\leq d)$ with the norm given by \eqref{space1}.
\begin{theorem}\label{thm2.20}
Let $(\mathcal{M},g)$ be a $d$-dimensional real hyperbolic manifold with $d \geqslant 2$.
Let $0<\delta<\min\left\{1,\, \dfrac{d}{8} \right\}$ and $0<\sigma<{\beta}= d-1+ \frac{\gamma_{d/\delta,d/\delta}+\gamma_{d/(2\delta),d/\delta}}{2}$. We suppose that $u(0)$ belongs to $L^p(\Gamma(T\mathcal{M}))\cap L^d(\Gamma(T\mathcal{M}))$ and $F$  belongs  to $ AAP(\r_+, (L^{\frac{dp}{d+\delta p}}\cap L^{\frac{d}{1+\delta}}\cap L^{\frac{d}{2\delta}})(\Gamma(T\mathcal{M}\otimes T\mathcal{M})))$ be  a forward asymptotically almost periodic function with respect to $t$ for $1<p\leq d$ and
$$P=\sup\limits_{t >0} \left( \| F(t)\|_{\frac{dp}{d+\delta p}} + \| F(t) \|_{\frac{d}{1+\delta}} + [h_d(t)]^{-\frac{1-\delta}{d}}e^{\sigma t}\|F(t)\|_{\frac{d}{2\delta}} \right)<+\infty.$$
Then, the following assertions hold
\begin{itemize}
\item[i)] If the norms $\|u(0)\|_p$, $\|u(0)\|_d$ and $P$ are sufficiently small, 
then equation \eqref{DivNavierStokes} has one and only one forward AAP- mild solution $\hat{u}$ on a small ball of  
$AAP(\r_+, (L^{p}\cap L^{d} \cap L^{\frac{d}{\delta}})(\Gamma(T\mathcal{M})))$.
\item[ii)] If $Q=\sup\limits_{t>0} e^{\sigma^* t} \left( \| F(t)\|_{\frac{dp}{d+\delta p}} + \| F(t) \|_{\frac{d}{1+\delta}} + [h_d(t)]^{-\frac{1-\delta}{d}} e^{\sigma t} \|F(t)\|_{\frac{d}{2\delta}} \right) <+\infty,$
then the above AAP- mild solution $\hat{u}(t)$ decays exponentially in the sense that for all $t>0$:
$$\left|\norm{\hat{u}(t)}\right| \lesssim e^{- \sigma^* t},$$
where $0<\sigma^*<\min\left\{ \beta_1,\, \hat{\beta}_1, \, \beta_2,\, \hat{\beta}_2 \right\}$ with
$\beta_2 = d-1+\gamma_{p,p}$, $\hat{\beta}_2 = d-1+\gamma_{d,d}$, $\beta_1 =d-1  + \frac{\gamma_{p,p}+\gamma_{dp/(1+\delta p),p}}{2}$, $\hat{\beta}_1 = d-1  + \frac{\gamma_{d,d} + \gamma_{d/(1+\delta),d}}{2}$ and $0<\sigma^*+\sigma< {\beta}=d-1+ \frac{\gamma_{d/\delta,d/\delta}+ \gamma_{d/(2\delta),d/\delta}}{2}$. Here, we denote
$$\left| \norm{\cdot}\right| = \norm{\cdot}_p + \norm{\cdot}_d + [h_d(t)]^{-\frac{1-\delta}{d}}e^{\sigma t}\norm{\cdot}_{\frac{d}{\delta}}.$$
\end{itemize}
Here, we use the notation $(L^r\cap L^s\cap L^z)(Z)=L^{r}(Z)\cap L^{s}(Z) \cap L^{z}(Z)$, where $Z$ is $\Gamma(T\mathcal{M})$ or $\Gamma(T\mathcal{M}\otimes T\mathcal{M})$.
\end{theorem} 
\begin{proof}
\noindent
$i)$ We consider the existence of AAP- mild solution of Equation \eqref{DivNavierStokes} on the ball centered at zero and radius $\rho$:
\begin{eqnarray}\label{bro}
\B_\rho^{AAP}&=&\left\{v\in AAP(\r_+, (L^p\cap L^d \cap L^{d/\delta})(\Gamma(T\mathcal{M}))):\right.\cr
&&\left. \norm{v}_{\mathcal{X}}=\sup_{t>0}\left( \|v(t)\|_{p} + \|v(t)\|_{d} + [h_d(t)]^{-\frac{1-\delta}{d}}e^{\sigma t}\norm{v(t)}_{\frac{d}{\delta}}\right) \le \rho \right\}.
\end{eqnarray}

For each $v\in \B_\rho^{AAP}$, we consider the linear equation 
\begin{equation}\label{ns1}
u(t)  = u(0) +\int_0^t e^{-(t-\tau)\A} \p\dive \left(-v(\tau)\otimes v(\tau) +F(\tau)\right) d\tau.
\end{equation} 
By applying Lemma \ref{Thm:linear}, Equation \eqref{ns1} has unique AAP- mild solution $u$ such that
\begin{eqnarray}\label{CoreEstimate}
\norm{u}_{\mathcal{X}} &\leqslant& \norm{u(0)}_{L^p\cap L^d} + M\norm{v}^2_{\mathcal{X}} + NP\cr
&\leqslant& \norm{u(0)}_{L^p\cap L^d} + M \rho^2 + NP\cr
&\leqslant& \rho,
\end{eqnarray}
if $\norm{u(0)}_{L^p\cap L^d},\, \rho$ and $P$ are small enough. Therefore, we can define a map $\Phi: \mathcal{X} = C_b(\r_+, Y) \to \mathcal{X}=C_b(\r_+, Y)$ as follows
\begin{equation}\label{defphi}
\begin{split}
\Phi(v)&=u
\end{split}
\end{equation}
if $\|u(0) \|_{L^p\cap L^d}$, $\rho$ and $P$ are small enough satisfying Estimate \eqref{CoreEstimate}. Hence, the map  $\Phi$ acts from  $\B_\rho^{AAP}$ into itself. Then, by Formula 
\eqref{mild:linear} with $-v\otimes v + F$ instead of $F$, the transformation  $\Phi$ is written by 
\begin{equation}\label{defphi1}
\Phi(v)(t) = u(0) + \int_0^t e^{-(t-\tau)\mathcal{A}} 
 \mathbb{P}\dive\left(-v(\tau) \otimes v (\tau) + F (\tau)\right) d\tau.
\end{equation}
Moreover,
for $v_1, v_2\in \B_\rho^{AAP}$, the function $u:=\Phi(v_1)-\Phi(v_2)$ becomes the unique forward AAP- mild solution to the equation 
$$\partial_tu  + \A u   = \p\dive (-v_1\otimes v_1+v_2\otimes v_2) = \p\dive(-v_1\otimes(v_1-v_2) - v_2\otimes(v_1-v_2)).$$
Thus, by \eqref{defphi1} and the same way to establish inequality \eqref{CoreEstimate} (using the same way as in the proof of Lemma \ref{Thm:linear}), we can prove that
\begin{eqnarray}\label{Core}
\norm{\Phi(v_1)-\Phi(v_2)}_{\mathcal{X}} &\leqslant 2M \rho \norm{v_1-v_2}_{\mathcal{X}}.
\end{eqnarray}
If $\rho$ is sufficiently small, then $2M\rho<1$ and $\Phi$ is a contraction on $\B_\rho^{AAP}$. Therefore, for these values of $\rho$ and $P$, there exists a unique fixed point $\hat{u}$ of $\Phi$, and by the definition of $\Phi$, this function $\hat{u}$ is an forward AAP- mild solution to Navier-Stokes equation \eqref{DivNavierStokes}. By using \eqref{Core}, the uniqueness of $\hat{u}$ in the small ball $\B_\rho^{AAP}$ is clearly.
 
We would like to note that, if we replace $\B_\rho^{AAP}$ by 
 $\B_\rho:=\{v\in \mathcal{X}=C_b(\r_+, Y) : \|v\|_{\mathcal{X}}\le \rho\}$, then  in a same way
as above, we can prove that for a bounded external force $F\in 
C_b(\r_+, X)$ (where $X= (L^{\frac{dp}{d+\delta p}}\cap L^{\frac{d}{1+\delta}}\cap L^{\frac{d}{2\delta}})(\Gamma(T\mathcal{M}\otimes T\mathcal{M})))$) satisfying 
$$P=\sup\limits_{t >0} \left( \| F(t)\|_{\frac{dp}{d+\delta p}} + \| F(t) \|_{\frac{d}{1+\delta}} + [h_d(t)]^{\frac{1-\delta}{d}}e^{\sigma t}\|F(t)\|_{\frac{d}{2\delta}} \right)$$
being small enough, there exists a unique bounded  solution $u\in \B_\rho$ with small radius $\rho$.
 
\medskip
$ii)$ In recent works \cite{HuyXuan2020,HuyXuan2021',XVQ2023}, the asymptotic behaviour of the mild solutions of incompressible Navier-Stoke equations \eqref{DivNavierStokes} on non-compact Einstein or generalized non-compact manifolds with negarive Ricci curvatures are obtained by using the cone inequality. Here, we give another proof for the exponential decay of $\hat{u}$ and $\nabla\hat{u}$ by using the Gronwall's inequality and the convergences of beta and gamma functions.

Following assertions $i)$ and $iii)$ Lemma \ref{estimates}, we have
\begin{eqnarray*}
\norm{\hat{u}(t)}_p &\leq& e^{-\beta_2 t}\norm{{u}(0)}_p + \int_0^t[h_d(t-\tau)]^{\frac{1+\delta}{d}}e^{-\beta_1(t-\tau)}\norm{\hat{u}(\tau)\otimes \hat{u}(\tau) + F(\tau)}_{\frac{dp}{d+\delta p}} d\tau \cr
&\leq& e^{-\beta_2 t}\norm{{u}(0)}_p + \int_0^t[h_d(t-\tau)]^{\frac{1+\delta}{d}}e^{-\beta_1(t-\tau)}\norm{\hat{u}(\tau)}_p \norm{\hat{u}(\tau)}_{\frac{d}{\delta}} d\tau \cr
&&+ \int_0^t[h_d(t-\tau)]^{\frac{1+\delta}{d}}e^{-\beta_1(t-\tau)}\norm{F(\tau)}_{\frac{dp}{d+\delta p}} d\tau\cr
&\leq& e^{-\beta_2 t}\norm{{u}(0)}_p + \rho\int_0^t[h_d(t-\tau)]^{\frac{1+\delta}{d}}[h_d(\tau)]^{\frac{1-\delta}{d}}e^{-\beta_1(t-\tau)}e^{-\sigma^*\tau}\norm{\hat{u}(\tau)}_p d\tau \cr
&&+ \int_0^t[h_d(t-\tau)]^{\frac{1+\delta}{d}}e^{-\beta_1(t-\tau)}\norm{F(\tau)}_{\frac{dp}{d+\delta p}} d\tau,
\end{eqnarray*}
where $\beta_2 = d-1+\gamma_{p,p}$ and $\beta_1 = d-1  + \frac{\gamma_{p,p}+\gamma_{dp/(d+\delta p),p}}{2}$.

Setting $y(\tau) = e^{\sigma^*\tau}\norm{\hat{u}(\tau)}_p$. Since $0<\sigma^*<\min\left\{\beta_1,\,\beta_2 \right\}$, we obtain that
\begin{eqnarray*}
y(t) &\leq& \norm{{u}(0)}_p + \rho\int_0^t [h_d(t-\tau)]^{\frac{1+\delta}{d}}[h_d(\tau)]^{\frac{1-\delta}{d}}e^{-(\beta_1-\sigma^*)(t-\tau)} y(\tau) d\tau \cr
&&+ \int_0^t[h_d(t-\tau)]^{\frac{1+\delta}{d}}e^{-(\beta_1-\sigma^*)(t-\tau)}e^{\sigma^* \tau}\norm{F(\tau)}_{\frac{dp}{d+\delta p}} d\tau \cr
&\leq&\norm{{u}(0)}_p + \rho\int_0^t [h_d(t-\tau)]^{\frac{1+\delta}{d}}[h_d(\tau)]^{\frac{1-\delta}{d}}e^{-(\beta_1-\sigma^*)(t-\tau)}y(\tau) d\tau \cr
&&+ \sup_{t>0} e^{\sigma^* t}\norm{F(t)}_{\frac{dp}{d+\delta p}}\int_0^t[h_d(t-\tau)]^{\frac{1+\delta}{d}}e^{-(\beta_1-\sigma^*)(t-\tau)}d\tau\cr
&\leq&\norm{{u}(0)}_p + \rho\int_0^t [h_d(t-\tau)]^{\frac{1+\delta}{d}}[h_d(\tau)]^{\frac{1-\delta}{d}}e^{-(\beta_1-\sigma^*)(t-\tau)}y(\tau) d\tau \cr
&&+ Q\int_0^t[h_d(t-\tau)]^{\frac{1+\delta}{d}}e^{-(\beta_1-\sigma^*)(t-\tau)}d\tau.
\end{eqnarray*}
The integrals
$$\int_0^t [h_d(t-\tau)]^{\frac{1+\delta}{d}}[h_d(\tau)]^{\frac{1-\delta}{d}}e^{-(\beta_1-\sigma^*)(t-\tau)} d\tau \leqslant M_5<+\infty,$$
$$\int_0^t [h_d(t-\tau)]^{\frac{1+\delta}{d}}e^{-(\beta_1-\sigma^*)(t-\tau)}d\tau \leqslant \tilde{M}_5<+\infty$$
converge by the mean of beta and gamma functions (see Appendix). Using the Gronwall's inequality we get 
$$|y(t)|\leqslant (\norm{u(0)}_p + Q\tilde{M}_5) e^{\rho M_5} \hbox{  for all  } t>0.$$
This leads to the exponentially asymptotic behaviour of $\hat{u}(t)$ in $L^p(\Gamma(T\mathcal{M}))$:
\begin{equation}\label{decay1}
\norm{\hat{u}(t)}_p \lesssim e^{-\sigma^* t} \hbox{  for all  } t>0.
\end{equation}
By the same way as above, we obtain also that
\begin{equation}\label{decay2}
\norm{\hat{u}(t)}_d \lesssim e^{-\sigma^* t} \hbox{  for all  } t>0.
\end{equation}

On the other hand, By using Assertion iii) in Lemma \ref{dispersive} and condition $0<\sigma^*+\sigma<{\beta}$, we have
\begin{eqnarray}\label{part311}
&&[h_d(t)]^{-\frac{1-\delta}{d}}e^{(\sigma^*+\sigma) t}\norm{\int_0^t e^{-(t-\tau)\A} \p\dive [-\hat{u}(\tau)\otimes \hat{u}(\tau)] d\tau}_{\frac{d}{\delta}} \cr
&\leqslant& C^{\frac{1-\delta}{d}}e^{(\sigma^*+\sigma) t}\int_0^t \norm{e^{-(t-\tau)\A} \p\dive (-\hat{u}(\tau)\otimes \hat{u}(\tau))}_{\frac{d}{\delta}} d\tau \hbox{   (because   } [h_d(t)]^{-\frac{1-\delta}{d}}\leqslant C^{\frac{1-\delta}{d}})\cr
&\leqslant& C^{\frac{1-\delta}{d}}e^{(\sigma^*+\sigma) t}\int_0^t [h_d(t-\tau)]^{\frac{1+\delta}{d}} e^{-{\beta}(t-\tau)} \|\hat{u}(\tau)\otimes \hat{u}(\tau)\|_{\frac{d}{2\delta}}d\tau\cr
&\leqslant& C^{\frac{1-\delta}{d}}e^{(\sigma^*+\sigma) t}\int_0^t [h_d(t-\tau)]^{\frac{1+\delta}{d}} e^{-{\beta}(t-\tau)} \|\hat{u}(\tau)\|_{\frac{d}{\delta}} \|\hat{u}(\tau)\|_{\frac{d}{\delta}}d\tau\cr
&\leqslant&  C^{\frac{1-\delta}{d}}\rho e^{\sigma^* t}\int_0^t [h_d(t-\tau)]^{\frac{1+\delta}{d}}[h_d(\tau)]^{\frac{1-\delta}{d}}e^{-({\beta}-\sigma)(t-\tau)}  \|\hat{u}(\tau)\|_{\frac{d}{\delta}} d\tau.
\end{eqnarray}
Moreover, we have
\begin{eqnarray}\label{part322}
&&[h_d(t)]^{-\frac{1-\delta}{d}}e^{(\sigma^*+\sigma) t}\norm{\int_0^t e^{-(t-\tau)\A} \p\dive F(\tau) d\tau}_{\frac{d}{\delta}} \cr
&\leqslant& C^{\frac{1-\delta}{d}}e^{(\sigma^*+\sigma) t}\int_0^t \norm{e^{-(t-\tau)\A} \p\dive F(\tau)}_{\frac{d}{\delta}}d\tau \hbox{   (because   } [h_d(t)]^{-\frac{1-\delta}{d}}\leqslant C^{\frac{1-\delta}{d}})\cr
&\leqslant& C^{\frac{1-\delta}{d}}e^{(\sigma^*+\sigma) t}\int_0^t [h_d(t-\tau)]^{\frac{1+\delta}{d}} e^{-{\beta}(t-\tau)} \|F(\tau)\|_{\frac{d}{2\delta}}d\tau\cr
&\leqslant& \sup_{t>0}[h_d(t)]^{\frac{1-\delta}{d}}e^{(\sigma^*+\sigma) t}\|F(t)\|_{\frac{d}{2\delta}}C^{\frac{1-\delta}{d}}\int_0^t [h_d(t-\tau)]^{\frac{1+\delta}{d}}[h_d(\tau)]^{\frac{1-\delta}{d}}e^{-({\beta}-\sigma^*-\sigma)(t-\tau)} d\tau \cr
&\leqslant& M_6Q,
\end{eqnarray}
where
$$C^{\frac{1-\delta}{d}}\int_0^t [h_d(t-\tau)]^{\frac{1+\delta}{d}}[h_d(\tau)]^{\frac{1-\delta}{d}}e^{-({\beta}-\sigma-\sigma^*)(t-\tau)} d\tau \leqslant M_6 <+\infty.$$
(this integral convergences by beta and gamma functions similar (see Appendix)).

Setting $z(t)= [h_d(t)]^{-\frac{1-\delta}{d}}e^{(\sigma^*+\sigma) t}\norm{\hat{u}(t)}_{\frac{d}{\delta}}$.
The condition $0<\sigma^*+\sigma< {\beta}<\tilde{\beta}_2 = d-1+ \gamma_{d,d/\delta}$ (for $0<\delta<\min\left\{ 1,\, \frac{d}{8}\right\}$) leads to 
$$[h_d(t)]^{-\frac{1-\delta}{d}}e^{(\sigma^*+\sigma) t}\norm{ e^{-t\mathcal{A}}u(0)}_{\frac{d}{\delta}} \leqslant \norm{u(0)}_d.$$
Combining this inequality with inequalities \eqref{part311} and \eqref{part322}, we get
\begin{equation}\label{Part333}
z(t) \leqslant (\norm{u(0)}_d + QM_6) + C^{\frac{1-\delta}{d}}\rho \int_0^t [h_d(t-\tau)]^{\frac{1+\delta}{d}}[h_d(\tau)]^{\frac{2(1-\delta)}{d}}e^{-({\beta}-\sigma^*-\sigma)(t-\tau)}  z(\tau) d\tau.
\end{equation}
The integral $\int_0^t [h_d(t-\tau)]^{\frac{1+\delta}{d}}[h_d(\tau)]^{\frac{2(1-\delta)}{d}}e^{-({\beta}-\sigma^*-\sigma)(t-\tau)} d\tau \leqslant \tilde{M}_6 <+\infty$, by beta and gamma functions (see Appendix). Using the Gronwall's inequality we get 
$$|z(t)|\leqslant (\norm{u(0)}_d + QM_6) e^{C^{\frac{1-\delta}{d}}\rho\tilde{M}_6} \hbox{  for all  } t>0.$$
This leads to
\begin{equation}\label{decay3}
[h_d(t)]^{-\frac{1-\delta}{d}}e^{\sigma t}\norm{\hat{u}(t)}_{\frac{d}{\delta}} \lesssim e^{-\sigma^* t} \hbox{  for all  } t>0.
\end{equation}
Since inequalities \eqref{decay1}, \eqref{decay2} and \eqref{decay3}, we obtain that
\begin{equation*}
\left|\norm{\hat{u}(t)}\right| \lesssim e^{-\sigma^*t}.
\end{equation*}

\end{proof}

\begin{remark}
\item [$\bullet$] As a direct consequence of the exponential decay we have the exponential stability of the small forward AAP- mild solution $\hat{u}$ of Equation \eqref{DivNavierStokes} on the real hyperbolic manifold $\mathcal{M}$ as follows: for another bounded mild solution $u\in \mathcal{X}$ of Equation \eqref{DivNavierStokes} if $\norm{u(0)-\hat{u}(0)}_{L^p\cap L^d}$ is small enough, then the following inequality holds
$$\left|\norm{\hat{u}(t)-u(t)}\right| \lesssim e^{-\sigma^* t}\norm{\hat{u}(0)-u(0)}_{L^p\cap L^d}$$
for all  $t>0$, where $0<\sigma^*<\min\left\{ \beta_1, \,\hat{\beta}_1,\, \beta_2, \hat{\beta}_2 \right\}$ and $0<\sigma^*+\sigma< {\beta}$. Here, we denote
$$\left|\norm{\cdot}\right| = \norm{\cdot}_{p} + \norm{\cdot}_{d} + [h_d(t)]^{-\frac{1-\delta}{d}}e^{\sigma t}\norm{\cdot}_{\frac{d}{\delta}}.$$
\item [$\bullet$] The results in this paper can be extended to the generalized non-compact Riemannian manifolds with negative Ricci curvatures such as Einstein and non-compact Riemannian manifolds which satisfy the conditions $(H_1)-(H_4)$ in \cite{Pi}. In these cases, the main obstacle comes from the non-commutation of the Hodge-Kodaira projection $\mathbb{P}$ with the Stokes semigroup $e^{-t\mathcal{A}}$, then the proofs of boundedness of mild solution $u(t)$ of the inhomogeneous Stokes equation and the preserved property of the solution operator $S$ on AAP- functions (as in the proofs of Lemma \ref{Thm:linear} and Theorem \ref{pest}) become more complicated. Our recent work \cite{HuyXuan2021'} on the existence and stability of periodic mild solutions of INSE on these manifolds can be useful to overcome this obstacle.\\
\item [$\bullet$] The results in this works and the ones in \cite{XVQ2023} answer completely the questions about the existence, uniqueness and asymptotic behaviour of forward AAP- solutions of Navier-Stokes equations on the real hyperbolic space $\mathbb{H}^d\,\, (d\geqslant 2)$ in the context of $L^p$-spaces for all $p>1$. However, in other spaces such as Euclidean space $\mathbb{R}^d$ (which is a flat space) or the non-compact Riemannian manifolds which have non-negative Ricci curvatures, these questions are remain open problems. The reason is that the Stokes semigroup $e^{-t\mathcal{A}}$ can be polynomial stable in $\mathbb{R}^d\,\, (d\geqslant 3)$ (or its unbounded domains) (see \cite[Lemma 2.1]{KoNa}) or even it can not stable in some non-compact manifolds with non-positive Ricci curvatures (see for example \cite{Ma}). This fact leads to the proof of the existence of bounded mild solution becomes more complicated and the Massera-type principle may no longer be true for all types of mild solution. For example of simple case, we consider the Stokes and Navier-Stokes equations on the flat space $\mathbb{R}^d\,\, (d\geqslant 3)$ (or its unbounded domains), there are some systematic works \cite{GHN,HiHuSe,Huy2014,Huy2018,HuyHaSacXuan2,HuyHaSacXuan2',HuyHaSacXuan3} which establish the boundedness of mild solutions of Stokes equation by using suitable Lorentz spaces $L^{p,\infty}(\mathbb{R}^d)$ instead of the $L^p$-spaces. Then, one can also prove the Massera-type principle for the existence of some types of mild solution (such as periodic, almost periodic and almost automorphic mild solutions) of Stokes and Navier-Stokes equations (see \cite{Huy2018,HuyHaSacXuan2,HuyHaSacXuan2',HuyHaSacXuan3}). However, since the semigroup $e^{-t\mathcal{A}}=e^{t\Delta}$ is polynomially stable, the solution operator $S$ may not preserve the forward asymptotically almost periodic property of the function $f$, even in Lorentz or other interpolation spaces. Therefore, the existence of a forward asymptotically almost periodic mild solution of the Stokes and Navier-Stokes equations on $\mathbb{R}^d$ (or an unbounded domain of $\mathbb{R}^d$) is still an open problem.

\end{remark}

\section{Appendix}\label{A}
In this part, we will prove the boundedness of the following integrals which are used in the previous sections
\begin{eqnarray*}
I_1 &=& \int_0^t [h_d(t-\tau)]^{\frac{1+\delta}{d}}[h_d(\tau)]^{\frac{1-\delta}{d}}e^{-\beta_1(t-\tau)}e^{-\sigma\tau}d\tau <+\infty,\cr
I_2 &=& \int_0^t [h_d(t-\tau)]^{\frac{1+\delta}{d}}[h_d(\tau)]^{\frac{1-\delta}{d}}e^{-\hat{\beta}_1(t-\tau)}e^{-\sigma\tau}d\tau <+\infty,\cr
I_3 &=& \int_0^t [h_d(t-\tau)]^{\frac{1+\delta}{d}}[h_d(\tau)]^{\frac{2(1-\delta)}{d}}e^{-({\beta}-\sigma)(t-\tau)}e^{-\sigma\tau}d\tau <+\infty,\cr
I_4 &=& \int_0^t [h_d(t-\tau)]^{\frac{1+\delta}{d}}[h_d(\tau)]^{\frac{1-\delta}{d}}e^{-({\beta} -\sigma)(t-\tau)} d\tau <+\infty,\cr
I_5 &=& \int_0^t [h_d(t-\tau)]^{\frac{1+\delta}{d}} [h_d(\tau)]^{\frac{1-\delta}{d}} e^{-(\beta_1-\sigma) (t-\tau)} d\tau <+\infty,\cr
\tilde{I}_5 &=& \int_0^t [h_d(t-\tau)]^{\frac{1+\delta}{d}}e^{-(\beta_1-\sigma)(t-\tau)}d\tau <+\infty,\cr
I_6 &=& \int_0^t [h_d(t-\tau)]^{\frac{1+\delta}{d}}[h_d(\tau)]^{\frac{1-\delta}{d}}e^{-({\beta}-\sigma-\sigma^*)(t-\tau)} d\tau <+\infty\cr
\tilde{I}_6 &=& \int_0^t [h_d(t-\tau)]^{\frac{1+\delta}{d}}[h_d(\tau)]^{\frac{2(1-\delta)}{d}}e^{-({\beta}-\sigma^*-\sigma)(t-\tau)} d\tau <+\infty,
\end{eqnarray*}
where t>0, $1<p\leqslant d$ and $0<\sigma<{\beta}$ for $I_3,\, I_4$; $0<\sigma<\beta_1$ for $I_5,\, \tilde{I}_5$; $0<\sigma+\sigma^*<{\beta}$ for $I_6,\, \tilde{I}_6$.

The boundedness of integrals $I_1,\, I_2,\, I_3,\, I_4$ are useful to the proof of Assertion i) in Lemma \ref{Thm:linear}. The boundedness of integrals $I_5,\, I_6, \, \tilde{I}_6$ serve to the proof of Assertion ii) in Theorem \ref{thm2.20}.
Moreover, if we consider the integral on $(-\infty,t]$ for all $t\in \mathbb{R}$, then we have the following bounded integrals 
\begin{eqnarray*}
\tilde{I}_1 &=& \int_{-\infty}^t [h_d(t-\tau)]^{\frac{1+\delta}{d}}\lambda^{-1}(\tau)e^{-\beta_1(t-\tau)}d\tau <+\infty,\cr
\tilde{I}_2 &=& \int_{-\infty}^t [h_d(t-\tau)]^{\frac{1+\delta}{d}}\lambda^{-1}(\tau)e^{-\hat{\beta}_1(t-\tau)}d\tau <+\infty,\cr
\tilde{I}_3 &=& \lambda(t)\int_{-\infty}^t [h_d(t-\tau)]^{\frac{1+\delta}{d}}\lambda^{-2}(\tau)e^{-{\beta}(t-\tau)}d\tau <+\infty,\cr
\tilde{I}_4 &=& \lambda(t)\int_{-\infty}^t [h_d(t-\tau)]^{\frac{1+\delta}{d}}\lambda^{-1}(\tau)e^{-{\beta}(t-\tau)} d\tau <+\infty.
\end{eqnarray*}
These boundedness are useful to prove the Assertion ii) in Lemma \ref{Thm:linear}.

The boundedness of integral $I_1,\, I_2,\, I_4,\, I_5$ and $I_6$ are similar. We prove only the boundedness of $I_5$. Indeed, we consider the following cases of $t$:

\underline{Case 1: $0<t<1$.} We have 
\begin{eqnarray*}
I_5 &=& \int_0^t [h_d(t-\tau)]^{\frac{1+\delta}{d}} [h_d(\tau)]^{\frac{1-\delta}{d}} e^{-(\beta_1-\sigma) (t-\tau)} d\tau\cr
&\leqslant&  \int_0^t (t-\tau)^{-\frac{1+\delta}{2}} \tau^{-\frac{1-\delta}{2}} d\tau\cr
&\leqslant& {\bf B}\left( \frac{1-\delta}{2},\,\frac{1+\delta}{2} \right)<+\infty,
\end{eqnarray*}
where ${\bf B}(\cdot,\, \cdot)$ is beta function.

\underline{Case 2: $1\leq t$.} We have 
\begin{eqnarray*}
I_5 &=& \int_0^t [h_d(t-\tau)]^{\frac{1+\delta}{d}} [h_d(\tau)]^{\frac{1-\delta}{d}} e^{-(\beta_1-\sigma) (t-\tau)} d\tau\cr
&=&  \int_0^1 [h_d(t-\tau)]^{\frac{1+\delta}{d}} [h_d(\tau)]^{\frac{1-\delta}{d}} e^{-(\beta_1-\sigma) (t-\tau)} d\tau \cr
&&+ \int_1^t [h_d(t-\tau)]^{\frac{1+\delta}{d}} [h_d(\tau)]^{\frac{1-\delta}{d}} e^{-(\beta_1-\sigma) (t-\tau)} d\tau\cr
&\leqslant&  \int_0^1 \left((t-\tau)^{-\frac{1+\delta}{2}}+1\right) \tau^{-\frac{1-\delta}{2}} e^{-(\beta_1 - \sigma) (t-\tau)} d\tau \cr
&&+ \int_1^t (t-\tau)^{-\frac{1+\delta}{2}} e^{-(\beta_1-\sigma) (t-\tau)} d\tau\cr
&\leqslant&  \int_0^t (t-\tau)^{-\frac{1+\delta}{2}}\tau^{-\frac{1-\delta}{2}} d\tau + \int_0^1 \tau^{-\frac{1-\delta}{2}} d\tau\cr
&&+ \int_0^t (t-\tau)^{-\frac{1+\delta}{2}} e^{-(\beta_1-\sigma) (t-\tau)} d\tau\cr
&\leqslant& {\bf B}\left( \frac{1-\delta}{2},\, \frac{1+\delta}{2} \right) + \frac{2}{1+\delta} + (\beta_1-\sigma)^{-\frac{1-\delta}{2}}\mathbf{\Gamma}\left(\frac{1-\delta}{2}\right)<+\infty,
\end{eqnarray*}
where $\mathbf{\Gamma}(\cdot)$ is gamma function.

The boundedness of integrals $I_3$ and $\tilde{I}_6$ are similar, we prove only for $I_3$ as follows:

\underline{Case 1: $0<t<1$.} We have 
\begin{eqnarray*}
I_3 &=& \int_0^t [h_d(t-\tau)]^{\frac{1+\delta}{d}} [h_d(\tau)]^{\frac{2(1-\delta)}{d}} e^{-({\beta}-\sigma) (t-\tau)} e^{-\sigma \tau} d\tau\cr
&\leqslant&  \int_0^t (t-\tau)^{-\frac{1+\delta}{2}} \tau^{-(1-\delta)} d\tau\cr
&\leqslant& t^{\frac{1-\delta}{2}}{\bf B}\left( \frac{1-\delta}{2},\,\delta \right)<+\infty,
\end{eqnarray*}
because $0<t<1$.

\underline{Case 2: $1\leq t$.} We have 
\begin{eqnarray*}
I_3 &=& \int_0^t [h_d(t-\tau)]^{\frac{1+\delta}{d}} [h_d(\tau)]^{\frac{2(1-\delta)}{d}} e^{-({\beta}-\sigma) (t-\tau)} e^{-\sigma \tau} d\tau\cr
&=&  \int_0^1 [h_d(t-\tau)]^{\frac{1+\delta}{d}} [h_d(\tau)]^{\frac{2(1-\delta)}{d}} e^{-({\beta}-\sigma) (t-\tau)} d\tau \cr
&&+ \int_1^t [h_d(t-\tau)]^{\frac{1+\delta}{d}} [h_d(\tau)]^{\frac{2(1-\delta)}{d}} e^{-({\beta}-\sigma) (t-\tau)} d\tau\cr
&\leqslant&  \int_0^1 \left((t-\tau)^{-\frac{1+\delta}{2}}+1\right) \tau^{-(1-\delta)}e^{-({\beta}-\sigma)\tau} d\tau \cr
&&+ \int_1^t (t-\tau)^{-\frac{1+\delta}{2}} e^{-({\beta}-\sigma) (t-\tau)} d\tau\cr
&\leqslant&  \left(2^{\frac{1+\delta}{2}}+1\right) \int_0^{1/2} \tau^{-(1-\delta)} d\tau + 2^{1-\delta}\int_{1/2}^1 \left((t-\tau)^{-\frac{1+\delta}{2}}+1\right) e^{-({\beta}-\sigma)\tau}d\tau\cr
&&+ \int_0^t (t-\tau)^{-\frac{1+\delta}{2}} e^{-({\beta}-\sigma) (t-\tau)} d\tau\cr
&\leqslant& \left(2^{\frac{1+\delta}{2}}+1\right)\frac{1}{\delta 2^\delta} + \frac{2^{1-\delta}}{{\beta}-\sigma}\left( e^{-\frac{{\beta}-\sigma}{2}}- e^{-({\beta}-\sigma)}\right) \cr
&&+ (2^{1-\delta}+1) ({\beta}-\sigma)^{-\frac{1-\delta}{2}}\mathbf{\Gamma}\left(\frac{1-\delta}{2}\right)<+\infty.
\end{eqnarray*}

The boundedness of integral $\tilde{I}_5$ is established as follows
\begin{eqnarray*}
\tilde{I}_5 &=& \int_0^t [h_d(t-\tau)]^{\frac{1+\delta}{d}}e^{-(\beta_1-\sigma)(t-\tau)}d\tau \cr
&\leqslant& \int_0^t \left( C^{\frac{1+\delta}{d}}(t-\tau)^{-\frac{1+\delta}{2}}+1\right)e^{-(\beta_1-\sigma)(t-\tau)}d\tau\cr
&\leqslant& C^{\frac{1+\delta}{d}} (\beta_1-\sigma)^{-\frac{1-\delta}{2}}{\bf \Gamma}\left(\frac{1-\delta}{2} \right) + \frac{1}{\beta_1-\sigma}\left( 1 - e^{-(\beta_1-\sigma)t}\right)<+\infty.
\end{eqnarray*}

Now, we consider the boundedness of integrals $\tilde{I}_1,\,\tilde{I}_2,\, \tilde{I}_3,\, \tilde{I}_4$ on $(-\infty,t]$ for all $t\in \mathbb{R}$. These integrals appear in the proof of Assertion ii) in Lemma \ref{Thm:linear}. The boundedness of $\tilde{I}_1,\, \tilde{I}_2,\, \tilde{I}_4$ are similar, we check for example the boundedness of $\tilde{I}_4$:
$$\tilde{I}_4 = \lambda(t)\int_{-\infty}^t [h_d(t-\tau)]^{\frac{1+\delta}{d}}\lambda^{-1}(\tau)e^{-{\beta}(t-\tau)}d\tau,$$
where we recall that for $0<\delta<\dfrac{\delta+1}{2}<\delta'<1$:
\begin{align*}
\lambda(t) = \begin{cases}
\displaystyle  [h_d(t)]^{-\frac{1-\delta}{d}}e^{\sigma t}&\hbox{   if  } t\geqslant 0, \cr
[h_d(|t|)]^{-\frac{1-\delta'}{d}} & \hbox{   if  } t<0.
\end{cases}
\end{align*}

\underline{Case 1: $t<0$.} We have
\begin{eqnarray*}
\tilde{I}_4 &\lesssim& \int_{-\infty}^t [h_d(t-\tau)]^{\frac{1+\delta}{d}}[h_d(|\tau|)]^{\frac{1-\delta'}{d}}e^{-{\beta}(t-\tau)}d\tau \hbox{    (because   } [h_d(|t|)]^{-\frac{1-\delta'}{d}}<C^{\frac{1-\delta'}{2}})\cr
&=&\int_{-\infty}^{t-1} [h_d(t-\tau)]^{\frac{1+\delta}{d}}[h_d(|\tau|)]^{\frac{1-\delta'}{d}}e^{-{\beta}(t-\tau)}d\tau\cr
&&+ \int_{t-1}^t [h_d(t-\tau)]^{\frac{1+\delta}{d}}[h_d(|\tau|)]^{\frac{1-\delta'}{d}}e^{-{\beta}(t-\tau)}d\tau\cr
&=&\int_{-\infty}^{t-1} e^{-{\beta}(t-\tau)}d\tau + \int_{t-1}^t (t-\tau)^{-\frac{1+\delta}{2}}[h_d(|\tau|)]^{\frac{1-\delta'}{d}}e^{-{\beta}(t-\tau)}d\tau\cr
&=& \frac{1}{{\beta}}e^{-{\beta}}+A,
\end{eqnarray*}
where
$$A= \int_{t-1}^t (t-\tau)^{-\frac{1+\delta}{2}}[h_d(|\tau|)]^{\frac{1-\delta'}{d}}e^{-{\beta}(t-\tau)}d\tau.$$
If $t\leqslant-1$, then
\begin{eqnarray*}
A&=&\int_{t-1}^t (t-\tau)^{-\frac{1+\delta}{2}}e^{-{\beta}(t-\tau)}d\tau \leqslant \frac{2}{1-\delta}<+\infty.
\end{eqnarray*}
If $-1<t<0$, then
\begin{eqnarray*}
A&=&\int_{t-1}^{-1} (t-\tau)^{-\frac{1+\delta}{2}}e^{-{\beta}(t-\tau)}d\tau\cr
&&+\int_{-1}^{t} (t-\tau)^{-\frac{1+\delta}{2}}|\tau|^{-\frac{1-\delta'}{2}}e^{-{\beta}(t-\tau)}d\tau\cr
&\leqslant& \frac{2}{1-\delta}\left( 1- (t+1)^{\frac{1-\delta}{2}}\right) + \int_{-1}^{t} (t-\tau)^{-\frac{1+\delta}{2}}|\tau|^{-\frac{1-\delta'}{2}}d\tau\cr
&=& \frac{2}{1-\delta}\left( 1- (t+1)^{\frac{1-\delta}{2}}\right) + \int_{0}^{t+1} \tau^{-\frac{1+\delta}{2}}|t-\tau|^{-\frac{1-\delta'}{2}}d\tau\cr
&=& \frac{2}{1-\delta}\left( 1- (t+1)^{\frac{1-\delta}{2}}\right) + \int_{0}^{t+1} \tau^{-\frac{1+\delta}{2}}(\tau-t)^{-\frac{1-\delta'}{2}}d\tau\cr
&\leqslant& \frac{2}{1-\delta}\left( 1- (t+1)^{\frac{1-\delta}{2}}\right) + \int_{0}^{t+1} \tau^{-\frac{1+\delta}{2}}\tau^{-\frac{1-\delta'}{2}}d\tau \hbox{   (because   }\tau+1>\tau-t>\tau>0)\cr
&=& \frac{2}{1-\delta}\left( 1- (t+1)^{\frac{1-\delta}{2}}\right) + \frac{2}{\delta'-\delta}(t+1)^{\frac{\delta'-\delta}{2}}<+\infty.
\end{eqnarray*}

\underline{Case 2: t>0.} We have 
\begin{eqnarray*}
\tilde{I}_4 &=& \int_{-\infty}^0 [h_d(t-\tau)]^{\frac{1+\delta}{d}}[h_d(|\tau|)]^{\frac{1-\delta'}{d}}e^{-{\beta}(t-\tau)}d\tau\cr
&&+ \int_0^t [h_d(t-\tau)]^{\frac{1+\delta}{d}}[h_d(\tau)]^{\frac{1-\delta'}{d}}e^{-({\beta}-\sigma)(t-\tau)}d\tau\cr
&=&B+C,
\end{eqnarray*}
where
$$B= \int_{-\infty}^0 [h_d(t-\tau)]^{\frac{1+\delta}{d}}[h_d(|\tau|)]^{\frac{1-\delta'}{d}}e^{-{\beta}(t-\tau)}d\tau,$$
$$C= \int_0^t [h_d(t-\tau)]^{\frac{1+\delta}{d}}[h_d(\tau)]^{\frac{1-\delta'}{d}}e^{-({\beta}-\sigma)(t-\tau)}d\tau.$$
The boundedness of $C$ is similar the boundedness of integrals $I_1,\,I_2,\,I_3,\,I_4$.
Now, we prove the boundedness of $B$ as follows
\begin{eqnarray*}
B &=& \int_{-\infty}^0 [h_d(t-\tau)]^{\frac{1+\delta}{d}}[h_d(|\tau|)]^{\frac{1-\delta'}{d}}e^{-{\beta}(t-\tau)}d\tau\cr
&=&\int_{-\infty}^{-1} [h_d(t-\tau)]^{\frac{1+\delta}{d}}e^{-{\beta}(t-\tau)}d\tau\cr
&&+ \int_{-1}^0 [h_d(t-\tau)]^{\frac{1+\delta}{d}}|\tau|^{\frac{1-\delta'}{d}}e^{-{\beta}(t-\tau)}d\tau\cr
&=&\int_1^\infty [h_d(\tau)]^{\frac{1+\delta}{d}}e^{-{\beta}\tau}d\tau + \int_{-1}^0 [h_d(t-\tau)]^{-\frac{1+\delta}{d}}|\tau|^{\frac{1-\delta'}{2}}e^{-{\beta}(t-\tau)}d\tau\cr
&\leqslant& \int_0^\infty \left(\tau^{-\frac{1+\delta}{2}}+1\right)e^{-{\beta}\tau}d\tau + \int_{-1}^0 [h_d(t-\tau)]^{-\frac{1+\delta}{d}}|\tau|^{\frac{1-\delta'}{2}}e^{-{\beta}(t-\tau)}d\tau\cr
&\leqslant& {\beta}^{-\frac{1-\delta}{2}}{\bf \Gamma}\left( \frac{1-\delta}{2}\right) + \frac{1}{{\beta}} + D,
\end{eqnarray*}
where
$$D= \int_{-1}^0 [h_d(t-\tau)]^{-\frac{1+\delta}{d}}|\tau|^{\frac{1-\delta'}{2}}e^{-{\beta}(t-\tau)}d\tau.$$
If $t>1$, we have
\begin{eqnarray*}
D&=& \int_{-1}^0 |\tau|^{-\frac{1-\delta'}{2}}e^{-{\beta}(t-\tau)}d\tau \leqslant \int_{-1}^0 |\tau|^{-\frac{1-\delta'}{2}}d\tau = \frac{2}{1+\delta}<+\infty.
\end{eqnarray*}
If $0<t\leqslant 1$, then
\begin{eqnarray*}
D&=& \int_{-1}^0 (t-\tau)^{-\frac{1+\delta}{2}}|\tau|^{-\frac{1-\delta'}{2}}e^{-{\beta}(t-\tau)}d\tau\cr
&\leqslant& \int_{-1}^0 (-\tau)^{-\frac{1+\delta}{2}}(-\tau)^{-\frac{1-\delta'}{2}} d\tau = \frac{2}{\delta'-\delta}<+\infty.
\end{eqnarray*}
The boundedness of $\tilde{I}_4$ is completed. The boundedness of $\tilde{I}_3$ is done by the same way as above, but we need the condition that $\dfrac{1+\delta}{2}<\delta'<1$ to guarantee that
\begin{equation*}
\int_{0}^{t+1}\tau^{-\frac{1+\delta}{2}}\tau^{-2\frac{1-\delta'}{2}}d\tau =  \int_{0}^{t+1}\tau^{\delta'-\frac{\delta+3}{2}} d\tau = \frac{2}{2\delta'-(1+\delta)}(t+1)^{\frac{2\delta'-(1+\delta)}{2}}<+\infty,
\end{equation*}
where   $-1<t<0$ and
\begin{equation*}
\int_{-1}^0 (-\tau)^{-\frac{1+\delta}{2}}(-\tau)^{-2\frac{1-\delta'}{2}}d\tau =  \int_{-1}^0 (-\tau)^{\delta'-\frac{\delta+3}{2}}d\tau = \frac{2}{2\delta'-(1+\delta)} <+\infty.
\end{equation*}
{\bf Acknowledgements}\\
This work is supported by Vietnam Institute for Advanced Study in Mathematics (VIASM) 2023.

\end{document}